# TRANSITION TIMES AND STOCHASTIC RESONANCE FOR MULTIDIMENSIONAL DIFFUSIONS WITH TIME PERIODIC DRIFT: A LARGE DEVIATIONS APPROACH

By Samuel Herrmann, Peter Imkeller and Dierk Peithmann

*Université Henri Poincaré Nancy I, Humboldt-Universität zu Berlin and Humboldt-Universität zu Berlin*

We consider potential type dynamical systems in finite dimensions with two meta-stable states. They are subject to two sources of perturbation: a slow external periodic perturbation of period $T$ and a small Gaussian random perturbation of intensity $\varepsilon$, and, therefore, are mathematically described as weakly time inhomogeneous diffusion processes. A system is in stochastic resonance, provided the small noisy perturbation is tuned in such a way that its random trajectories follow the exterior periodic motion in an optimal fashion, that is, for some optimal intensity $\varepsilon(T)$. The physicists' favorite, measures of quality of periodic tuning—and thus stochastic resonance—such as spectral power amplification or signal-to-noise ratio, have proven to be defective. They are not robust w.r.t. effective model reduction, that is, for the passage to a simplified finite state Markov chain model reducing the dynamics to a pure jumping between the meta-stable states of the original system. An entirely probabilistic notion of stochastic resonance based on the transition dynamics between the domains of attraction of the meta-stable states—and thus failing to suffer from this robustness defect—was proposed before in the context of one-dimensional diffusions. It is investigated for higher-dimensional systems here, by using extensions and refinements of the Freidlin–Wentzell theory of large deviations for time homogeneous diffusions. Large deviations principles developed for weakly time inhomogeneous diffusions prove to be key tools for a treatment of the problem of diffusion exit from a domain and thus for the approach of stochastic resonance via transition probabilities between meta-stable sets.

**Introduction.** The ubiquitous phenomenon of stochastic resonance has been studied by physicists for about 20 years and recently discovered in numerous areas of natural sciences. Its investigation took its origin in a toy model from climatology, which may serve to explain some of its main features.

To give a qualitative explanation for the almost periodic recurrence of cold and warm ages (glacial cycles) in paleoclimatic data, Nicolis [12] and Benzi, Sutera and Vulpiani [3] proposed a simple stochastic climate model based on an energy balance equation for the averaged global temperature $T(t)$ at time $t$. The balance between averaged absorbed and emitted radiative energies leads to a deterministic differential equation for $T(t)$ of the form

$$(0.1) \qquad \dot{T}(t) = b(Q(t), T(t)).$$

The solar constant $Q(t)$ fluctuates periodically at a very low frequency of $10^{-5}$ times per year due to periodic changes of the Earth orbit's eccentricity (Milankovich cycles), which coincide roughly with the observed frequency of ice and warm ages. Under reasonable assumptions, for frozen $q$, the nonlinear function $b(q, T)$ describes the force associated with a double well potential possessing two stable temperature states which represent cold and warm ages. As $Q$ varies periodically, these states become meta-stable and are moved periodically by $Q(t)$. Most importantly, *transitions between these states are impossible*. Only the addition of a *stochastic forcing* allows for spontaneous transitions between the meta-stable climate states, thus explaining roughly transition mechanisms leading to glacial cycles.

In general, trajectories of the solutions of differential equations of this type, subject to two independent sources of perturbation, an exterior periodic one of period $T$, and a random one of intensity $\varepsilon$, say, will exhibit some kind of randomly periodic behavior, reacting to the periodic input forcing and eventually amplifying it. The problem of *optimal tuning* at large periods $T$ consists in finding a noise amplitude $\varepsilon(T)$ (the *resonance point*) which supports this amplification effect in a *best possible way*. During the last 20 years, various concepts of measuring the quality of periodic tuning to provide a criterion for optimality have been discussed and proposed in many applications from a variety of branches of natural sciences (see [7] for an overview). Its mathematical treatment started only very recently, and criteria for finding an optimal tuning are still under discussion.

The first approach toward a mathematically precise understanding of stochastic resonance was done by Freidlin [5]. Using large deviations theory, he explains basic periodicity properties of the trajectories in the large period (small noise) limit by the effect of deterministic quasi-periodic motion, but fails to account for optimal tuning. The most prominent quality measures for periodic tuning from the physics literature, the *signal-to-noise ratio* and the *spectral power amplification coefficient* (SPA), were investigated in a mathematically precise way in Pavlyukevich's thesis [13], and seen to have a serious drawback. Due to the high complexity of original systems, when calculating the optimal noise intensity, physicists usually pass to the effective dynamics of some kind of simple caricature of the system reducing the diffusion dynamics to the pure *inter well motion* (see, e.g., [11]). The reduced dynamics are represented by a continuous time two state Markov chain. Surprisingly, due to the importance of small *intra well fluctuations*, the tuning and resonance pattern of the Markov chain model may differ essentially from the resonance picture of the diffusion. It was this lack of robustness against model reduction which motivated Herrmann and Imkeller [9] to look for different measures of quality of periodic tuning for diffusion trajectories, retaining only the rough inter well motion of the diffusion. The measure they treat in the setting of *one-dimensional diffusion processes* subject to periodic forcing of small frequency is related to the transition probability during a fixed time window of exponential length, the position of which is tracked by a parameter of period length in which maximization is performed to account for optimal tuning.

The subject of the present paper is to continue our previous work in the general setting of *finite-dimensional diffusion processes*. Our approach of stochastic resonance thereby is based on the same robust probabilistic notion of periodic tuning. This extension is by no means obvious, since the multidimensional problem requires entirely new methods. We recall at this point that in [9] methods of investigation of stochastic tuning were heavily based on comparison arguments which are not an appropriate tool from dimension 2 on. Time inhomogeneous diffusion processes, such as the ones under consideration, were compared to piecewise homogeneous diffusion processes by freezing the potential's time dependence on small intervals. We study a dynamical system in $d$-dimensional Euclidean space perturbed by a $d$-dimensional Brownian motion $W$, that is, we consider the solution of the stochastic differential equation

$$(0.2) \qquad dX_t^\varepsilon = b\left(\frac{t}{T}, X_t^\varepsilon\right) dt + \sqrt{\varepsilon}\, dW_t, \qquad t \geq 0.$$

One of the system's important features is that its inhomogeneity is weak in the sense that the drift depends on time only through a re-scaling by the time parameter $T = T^\varepsilon$ which will be assumed to be exponentially large in $\varepsilon$. This corresponds to the situation in [9] and is motivated by the well-known Kramers–Eyring law which was mathematically underpinned by the Freidlin–Wentzell theory of large deviations [6]. The law roughly states that the expected time it takes for a homogeneous diffusion to leave a local attractor, for example, across a potential wall of height $\frac{v}{2}$, is given to exponential order by $T^\varepsilon = \exp(\frac{v}{\varepsilon})$. Hence, only in exponentially large scales of the form $T^\varepsilon = \exp(\frac{\mu}{\varepsilon})$ we can expect to see effects of transitions between different domains of attraction. $b$ is assumed to be one-periodic w.r.t. time. The deterministic system $\dot\xi_t = b(s, \xi_t)$ with *frozen* time parameter $s$ is supposed to have two domains of attraction that do not depend on $s \geq 0$. In the "classical" case of a drift derived from a potential, $b(t,x) = -\nabla_x U(t,x)$ for some potential function $U$, equation (0.2) describes the motion of a Brownian particle in a $d$-dimensional time inhomogeneous double-well potential.

Since our stochastic resonance criterion is based on transition times between the two meta-stable sets of the system, our analysis relies on a suitable notion of transition or exit time. The Kramers–Eyring formula suggests to consider the parameter $\mu$ from $T^\varepsilon = \exp(\frac{\mu}{\varepsilon})$ as a natural measure of scale. Therefore, if at time $s$ the system needs energy $e(s)$ to leave some meta-stable set, an exit from that set should occur at time

$$a_\mu = \inf\{t \geq 0 : e(t) \leq \mu\}$$

in the diffusion's natural time scale. If $a_\mu^i$ are the transition times for the two domains of attraction numbered $i = \pm 1$, we look at the probabilities of transitions between them within a time window $[(a_\mu^i - h)T^\varepsilon, (a_\mu^i + h)T^\varepsilon]$ for small $h > 0$. Assume for this purpose that the two corresponding meta-stable points are given by $x_i, i = \pm 1$, and denote by $\tau_\varrho^{-i}$ the random time at which the diffusion reaches the $\varrho$-neighborhood $B_\varrho(x_{-i})$ of $x_{-i}$. Then we use the following quantity to measure the quality of periodic tuning:

$$\mathcal{M}(\varepsilon, \mu) = \min_{i=\pm 1} \sup_{x \in B_\varrho(x_i)} \mathbb{P}_x(\tau_\varrho^{-i} \in [(a_\mu^i - h)T, (a_\mu^i + h)T]),$$

the minimum being taken in order to account for transitions back and forth. In order to exclude trivial or chaotic transition behavior, the scale parameter $\mu$ has to be restricted to an interval $I_R$ of reasonable values which we call *resonance interval*. With this measure of quality, the stochastic resonance point may be determined as follows. We fix the window width parameter $h > 0$ and maximize $\mathcal{M}(\varepsilon, \mu)$ in $\mu$ asymptotically as $\varepsilon \to 0$, or, equivalently, as $T^\varepsilon = \exp \frac{\mu}{\varepsilon} \to \infty$. If the asymptotic maximizer exists and is reached for the time scale $\mu_0(h)$, we call the eventually existing limit $\mu_0 = \lim_{h \to 0} \mu_0(h)$ resonance point. This is the optimal tuning w.r.t. our quality measure, that is, the best asymptotic relation between the noise amplitude $\varepsilon$ and the period $T^\varepsilon$.

To calculate $\mu_0(h)$ for fixed positive $h$, we use large deviations techniques. In fact, our main result contains a formula which states that

$$\lim_{\varepsilon \to 0} \varepsilon \log\{1 - \mathcal{M}(\varepsilon, \mu)\} = \max_{i=\pm 1}\{\mu - e_i(a_\mu^i - h)\}.$$

We show that this asymptotic relation holds uniformly w.r.t. $\mu$ on compact subsets of $I_R$, a fact which enables us to perform a maximization and find $\mu_0(h)$. The techniques needed to prove our main result feature extensions and refinements of the fundamental large deviations theory for time homogeneous diffusions by Freidlin–Wentzell [6]. We prove a large deviations principle for the inhomogeneous diffusion (0.2) and strengthen this result to get uniformity in system parameters. Similarly to the time homogeneous case, where large deviations theory is applied to the problem of *diffusion exit* culminating in a mathematically rigorous proof of the Kramers–Eyring law, we study the problem of diffusion exit from a domain which is carefully chosen in order to allow for a detailed analysis of transition times. The main idea behind our analysis is that the natural time scale is so large that rescaling in these units essentially leads to an asymptotic freezing of the time inhomogeneity, which has to be carefully studied, to hook up to the theory of large deviations of time homogeneous diffusions.

The material in the paper is organized as follows. Section 1 is devoted to the careful extension of large deviations theory to diffusions with slow time inhomogeneity. The most useful result for the subsequent analysis of exit times is Proposition 1.8, with a large deviations principle for slowly time dependent diffusions, uniform with respect to a system parameter. In Section 2 upper and lower bounds for the asymptotic exponential exit rate from domains of attraction for slowly time dependent diffusions are derived. The main result Theorem 2.3 combines them. Section 3 is concerned with developing the resonance criterion and computing the resonance point from the results of the preceding section.

**1. Large deviations for diffusion processes.** Let us now consider dynamical systems driven by slowly time dependent vector fields, perturbed by Gaussian noise of small intensity. We shall be interested in their large deviation behavior. Due to the slow time inhomogeneity, the task we face is not covered by the classical theory presented in [4] and [6]. For this reason, we shall have to extend the theory of large deviations for randomly perturbed dynamical systems developed by Freidlin and Wentzell [6] to drift terms depending in a weak form to be made precise below on the time parameter. Before doing so in the second subsection, we shall recall the classical results on time homogeneous diffusions in the following brief overview.

1.1. *The time homogeneous case*: *classical results.* For a more detailed account of the following well-known theory, see [4] or [6].

We consider the family of $\mathbb{R}^d$-valued processes $X^\varepsilon$, $\varepsilon > 0$, defined by

$$(1.1) \qquad dX_t^\varepsilon = b(X_t^\varepsilon)\,dt + \sqrt{\varepsilon}\,dW_t, \qquad X_0^\varepsilon = x_0 \in \mathbb{R}^d,$$

on a fixed time interval $[0,T]$, where $b$ is Lipschitz continuous and $W$ is a $d$-dimensional Brownian motion. This family of diffusion processes satisfies in the small noise limit, that is, as $\varepsilon \to 0$, a large deviations principle (LDP) in the space $C_{0T} := C([0,T], \mathbb{R}^d)$ equipped with the topology of uniform convergence induced by the metric $\rho_{0T}(\varphi, \psi) := \sup_{0 \leq t \leq T} \|\varphi_t - \psi_t\|$, $\varphi, \psi \in C_{0T}$. The *rate function* or *action functional* is given by $\bar{I}_{0T}^{x_0} : C_{0T} \to [0, +\infty]$,

$$(1.2) \quad I_{0T}^{x_0}(\varphi) = \begin{cases} \frac{1}{2} \int_0^T \|\dot{\varphi}_t - b(\varphi_t)\|^2 \, dt, \\ \qquad \text{if } \varphi \text{ is absolutely continuous and } \varphi_0 = x_0, \\ +\infty, \quad \text{otherwise.} \end{cases}$$

Moreover, $I_{0T}^{x_0}$ is a *good* rate function, that is, it has compact level sets. The LDP for this family of processes is mainly obtained as an application of the contraction principle to the LDP for the processes $\sqrt{\varepsilon} W$, $\varepsilon > 0$. More precisely, in the language of Freidlin and Wentzell, the functional $I_{0T}^{x_0}$ is the *normalized* action functional corresponding to the normalizing coefficient $\frac{1}{\varepsilon}$. In the sequel we will not consider scalings other than this one. We have $I_{0T}^{x_0}(\varphi) < \infty$ if and only if $\varphi$ belongs to the Cameron–Martin space of absolutely continuous functions with square integrable derivatives starting at $x_0$, that is,

$$\varphi \in H_{0T}^{x_0} := \left\{ f : [0,T] \to \mathbb{R}^d \Big| f(t) = x_0 + \int_0^t g(s) \, ds \text{ for some } g \in L^2([0,T]) \right\}.$$

We omit the superscript $x_0$ whenever there is no confusion about the initial condition we are referring to.

Observe that $I_{0T}(\varphi) = 0$ means that $\varphi$ (up to time $T$) is a solution of the deterministic equation

$$(1.3) \quad \dot{\xi} = b(\xi),$$

so $I_{0T}(\varphi)$ is essentially the $L^2$-deviation of $\varphi$ from the deterministic solution $\xi$. The *cost function* $V$ of $X^\varepsilon$, defined by

$$V(x, y, t) = \inf\{I_{0t}(\varphi) : \varphi \in C_{0t}, \varphi_0 = x, \varphi_t = y\},$$

takes into account all continuous paths connecting $x, y \in \mathbb{R}^d$ in a fixed time interval of length $t$, and the *quasi-potential*

$$V(x, y) = \inf_{t > 0} V(x, y, t)$$

describes the cost of $X^\varepsilon$ going from $x$ to $y$ eventually. In the potential case, $V$ agrees up to a constant with the potential energy to spend in order to pass from $x$ to $y$ in the potential landscape, hence, the term quasi-potential.

The classical LDP due to Freidlin and Wentzell requires the usual global Lipschitz and linear growth conditions from the standard existence and uniqueness results for SDE. In our setting the coefficients will (in general) not be globally Lipschitz since the drift is given by a potential gradient. An extension to locally Lipschitz and $\varepsilon$-dependent drift terms was provided by Azencott [1]. The following proposition is a special case of [1], Chapter III, Theorem 2.13. See also [2], Theorem 2.1.

PROPOSITION 1.1. *Assume that the equation* (1.1) *has a unique strong solution that never explodes and that the drift is locally Lipschitz. Then $X^\varepsilon$ satisfies on any time interval $[0,T]$ a large deviations principle with good rate function $I_{0T}$. Furthermore, the LDP for $X^\varepsilon$ holds uniformly w.r.t. the initial condition of the diffusion. More precisely, if $\mathbb{P}_y(X^\varepsilon \in \cdot)$ denotes the law of the diffusion $X^\varepsilon$ starting in $y \in \mathbb{R}^d$ and $K \subset \mathbb{R}^d$ is compact, we have, for any closed $F \subset C_{0T}$,*

$$(1.4) \qquad \limsup_{\varepsilon \to 0} \varepsilon \log \sup_{y \in K} \mathbb{P}_y(X^\varepsilon \in F) \leq -\inf_{y \in K} \inf_{\varphi \in F} I_{0T}^y(\varphi)$$

*and for any open $G \subset C_{0T}$,*

$$(1.5) \qquad \liminf_{\varepsilon \to 0} \varepsilon \log \inf_{y \in K} \mathbb{P}_y(X^\varepsilon \in G) \geq -\sup_{y \in K} \inf_{\varphi \in G} I_{0T}^y(\varphi).$$

REMARK 1.2. (i) A sufficient condition for the existence of a nonexploding and unique strong solution is a locally Lipschitz drift term $b$ which satisfies

$$(1.6) \qquad \langle x, b(x) \rangle \leq \gamma(1 + \|x\|^2) \qquad \text{for all } x \in \mathbb{R}^d$$

for some constant $\gamma > 0$ (see [17], Theorem 10.2.2). This still rather weak condition is obviously satisfied if $\langle x, b(x) \rangle \leq 0$ for large enough $x$, which means that $b$ contains a component that pulls $X$ back to the origin.

(ii) A strengthening of condition (1.6) ensuring superlinear growth will be used in subsequent sections. In that case, the laws of $(X^\varepsilon)$ are exponentially tight, and $I_{0T}$ is a good rate function. Recall that the laws of $(X^\varepsilon)$ are *exponentially tight* if there exist some $R_0 > 0$ and a positive function $\varphi$ satisfying $\lim_{x \to \infty} \varphi(x) = +\infty$ such that

$$(1.7) \qquad \limsup_{\varepsilon \to 0} \varepsilon \log \mathbb{P}(\sigma_R^\varepsilon \leq T) \leq -\varphi(R) \qquad \text{for all } R \geq R_0.$$

Here $\sigma_R^\varepsilon$ denotes the first time that $X^\varepsilon$ exits from $B_R(0)$.

1.2. *General results on weakly time inhomogeneous diffusions.* Let us now come to inhomogeneous diffusions with slowly time dependent drift coefficients. For our understanding of stochastic resonance effects of dynamical systems with slow time dependence, we have to adopt the large deviations results of the previous subsection to diffusions moving in potential landscapes with different valleys slowly and periodically changing their depths and positions. In this subsection we shall extend the large deviations results of Freidlin and Wentzell to time inhomogeneous diffusions which are almost homogeneous in the small noise limit, so that, in fact, we are able to compare to the large deviation principle for time homogeneous diffusions. The result we present in this subsection is not strong enough for the treatment

8    S. HERRMANN, P. IMKELLER AND D. PEITHMANNof stochastic resonance (one needs uniformity in some of the system parameters), but it most clearly exhibits the idea of the approach, which is why we state it here. Consider the family $X^\varepsilon$, $\varepsilon > 0$, of solutions of the SDE

$$(1.8) \qquad dX_t^\varepsilon = b^\varepsilon(t, X_t^\varepsilon)\,dt + \sqrt{\varepsilon}\,dW_t, \qquad t \geq 0,\ X_0^\varepsilon = x_0 \in \mathbb{R}^d.$$

We assume that (1.8) has a global strong solution for all $\varepsilon > 0$. Our main large deviations result for diffusions for which time inhomogeneity fades out in the small noise limit is summarized in the following proposition. The $\varepsilon$-dependence of the drift term was assumed in the same way in [1], Chapter III, Theorem 2.13 and [2], Theorem 2.1. See also [14].

PROPOSITION 1.3 (Large deviations principle). *Assume that the drift of the SDE* (1.8) *satisfies*

$$(1.9) \qquad \lim_{\varepsilon \to 0} b^\varepsilon(t, x) = b(x)$$

*for all $t \geq 0$, uniformly w.r.t. $x$ on compact subsets of $\mathbb{R}^d$, for some locally Lipschitz function $b : \mathbb{R}^d \to \mathbb{R}^d$. Assume that the time homogeneous diffusion $Y^\varepsilon$ associated to the limiting drift $b$ [i.e., the solution of (1.1) with the same initial condition $x_0$] does not explode.*

*Then $(X^\varepsilon)$ satisfies a large deviations principle on any finite time interval $[0, T]$ with good rate function $I_{0T}$ given by (1.2).*

PROOF. For notational convenience, we drop the $\varepsilon$-dependence of $X$ and $Y$. We shall prove that $X$ and $Y$ are exponentially equivalent, that is, for any $\delta > 0$, we have

$$(1.10) \qquad \limsup_{\varepsilon \to 0} \varepsilon \log \mathbb{P}(\rho_{0T}(X, Y) > \delta) = -\infty.$$

In order to verify this, fix some $\delta > 0$, and observe that

$$\|X_t - Y_t\| \leq \int_0^t \|b^\varepsilon(u, X_u) - b(X_u)\|\,du + \int_0^t \|b(X_u) - b(Y_u)\|\,du.$$

For $R > 0$, let $\tau_R := \inf\{t \geq 0 : X_t \notin B_R(x_0)\}$, let $\tilde{\tau}_R$ be defined similarly with $X$ replaced by $Y$, and $\sigma_R := \tau_R \wedge \tilde{\tau}_R$. The local Lipschitz continuity of $b$ implies the existence of some constant $K_R(x_0)$ such that $\|b(x) - b(y)\| \leq K_R(x_0)\|x - y\|$ for $x, y \in B_R(x_0)$. An application of Gronwall's lemma yields

$$\rho_{0T}(X, Y) \leq e^{K_R(x_0)T} \int_0^T \|b^\varepsilon(u, X_u) - b(X_u)\|\,du \qquad \text{on } \{\sigma_R > T\}.$$

Due to uniform convergence, for any $\eta > 0$, we can find some $\varepsilon_0 > 0$ s.t.

$$\sup_{x \in B_R(x_0)} \|b^\varepsilon(t, x) - b(x)\| \leq \eta \qquad \text{for } t \in [0, T], \varepsilon < \varepsilon_0.$$

This implies

(1.11) $$\rho_{0T}(X,Y) \leq \eta T e^{K_R(x_0)T} \qquad \text{for } \varepsilon < \varepsilon_0 \text{ on } \{\sigma_R > T\}.$$

By choosing $\eta$ small enough s.t. $\rho_{0T}(X,Y) \leq \delta/2$ on $\{\sigma_R > T\}$ [i.e., $X$ and $Y$ are very close before they exit from $B_R(x_0)$], we see that, for $\varepsilon < \varepsilon_0$,

$$\mathbb{P}(\rho_{0T}(X,Y) > \delta) \leq \mathbb{P}(\tau_R \leq T) + \mathbb{P}(\tilde{\tau}_R \leq T).$$

Since $X$ and $Y$ are close within the ball $B_R(x_0)$, we deduce that if $X$ escapes from $B_R(x_0)$ before time $T$, then $Y$ must at least escape from $B_{R/2}(x_0)$ before time $T$ (if $R > \delta$). So we have $\mathbb{P}(\rho_{0T}(X,Y) > \delta) \leq \mathbb{P}(\tilde{\tau}_{R/2} \leq T)$ for $\varepsilon < \varepsilon_0$. Hence, the LDP for $Y$ gives

$$\limsup_{\varepsilon \to 0} \varepsilon \log \mathbb{P}(\rho_{0T}(X,Y) > \delta) \leq -\inf_{0 \leq t \leq T, \|y - x_0\| \geq R/2} V(x_0, y, t).$$

Sending $R \to \infty$ yields the desired result (see Theorem 4.2.13 in [4]). □

It is easy to see that the uniformity w.r.t. the diffusion's initial condition also holds for the weakly inhomogeneous process $X^\varepsilon$ of this proposition. One only has to carry over Proposition 5.6.14 in [4], which is easily done using some Gronwall argument. Then the proof of the uniformity is the same as in the homogeneous case (see [4], Corollary 5.6.15). We omit the details.

1.3. *Weak inhomogeneity through slow periodic variation.* In this subsection we shall deal with some particular diffusions for which the drift term is subject to a very slow periodic time inhomogeneity. More precisely, we shall be concerned with solutions of the following stochastic differential equation taking their values in $d$-dimensional Euclidean space, driven by a $d$-dimensional Brownian motion $W$ of intensity $\varepsilon$:

(1.12) $$dX_t^\varepsilon = b\left(\frac{t}{T^\varepsilon}, X_t^\varepsilon\right) dt + \sqrt{\varepsilon}\, dW_t, \qquad t \geq 0, X_0 = x_0 \in \mathbb{R}^d.$$

Here $T^\varepsilon$ is a time scale parameter which tends to infinity as $\varepsilon \to 0$. In the subsequent sections, we shall assume that $T^\varepsilon$ is exponentially large, in fact,

(1.13) $$T^\varepsilon = \exp\frac{\mu}{\varepsilon} \qquad \text{with } \mu > 0.$$

The drift $b(t,x)$ of (1.12) is a time-periodic function of period one. Concerning its regularity properties, we suppose it to be locally Lipschitz in both variables, that is, for $R > 0$, $x \in \mathbb{R}^d$, there are constants $K_R(x)$ and $\kappa_R(x)$ such that

(1.14) $$\|b(t,y_1) - b(t,y_2)\| \leq K_R(x)\|y_1 - y_2\|,$$
(1.15) $$\|b(t,y) - b(s,y)\| \leq \kappa_R(x)|t - s|$$

for all $y, y_1, y_2 \in B_R(x)$ and $s, t \geq 0$. Furthermore, we shall assume that the drift term forces the diffusion to stay in compact sets for long times in order to get sufficiently "small" level sets. We suppose that there are constants $\eta$, $R_0 > 0$ such that

$$\langle x, b(t, x) \rangle < -\eta \|x\| \tag{1.16}$$

for $t \geq 0$ and $\|x\| \geq R_0$. This condition is stronger than (1.6), so the existence of a unique strong and nonexploding solution is again guaranteed. Moreover, this growth condition implies the exponential tightness of the diffusion (see Proposition 1.4 for the precise asymptotics).

1.3.1. *Boundedness of the diffusion.* The aim of this subsection is to exploit the consequences of the growth condition (1.16). In fact, it implies that the diffusion (1.12) cannot leave compact sets in the small noise limit. For positive $\varepsilon$, it stays for a long time in bounded domains. In the following proposition we shall make precise how the law of the exit time from bounded domains depends on $\varepsilon$. The arguments are borrowed from the framework of self-attracting diffusions; see [15] or [10].

For $R > 0$, let $\sigma_R^\varepsilon := \inf\{t \geq 0 : \|X_t^\varepsilon\| \geq R\}$ denote the first exit time from the ball $B_R(0)$.

PROPOSITION 1.4. *Let $\delta > 0$, and let $r : (0, \delta) \to (0, \infty)$ be a function satisfying $\lim_{\varepsilon \to 0} \frac{\varepsilon}{r(\varepsilon)} = 0$. There exist $R_1, \varepsilon_1 > 0$ and $C > 0$ such that, for $R \geq R_1$, $\varepsilon < \varepsilon_1$,*

$$\mathbb{P}_x(\sigma_R^\varepsilon \leq r(\varepsilon)) \leq C\eta^2 \frac{r(\varepsilon)}{\varepsilon} e^{-\eta R/\varepsilon} \qquad \text{for } \|x\| \leq \frac{R}{2}. \tag{1.17}$$

PROOF. For convenience of notation, we suppress the superscript $\varepsilon$ in $X^\varepsilon, \sigma_R^\varepsilon$ and so on. Choose a $C^2$-function $h : \mathbb{R}^d \to \mathbb{R}$ s.t. $h(x) = \|x\|$ for $\|x\| \geq R_0$ and $h(x) \leq R_0$ for $\|x\| \leq R_0$, where $R_0$ is the constant given in the growth condition (1.16). By Itô's formula, we have

$$h(X_t) = h(x) + \sqrt{\varepsilon} \int_0^t \nabla h(X_s) \, dW_s$$
$$+ \int_0^t \left\langle \nabla h, b\left(\frac{s}{T^\varepsilon}, \cdot\right) \right\rangle (X_s) \, ds + \frac{\varepsilon}{2} \int_0^t \triangle h(X_s) \, ds.$$

Let $\xi_t := \int_0^t \|\nabla h(X_s)\|^2 \, ds$, that is, $\xi_t$ is the quadratic variation of the continuous local martingale $M_t := \int_0^t \nabla h(X_s) \, dW_s, t \geq 0$. Since $\nabla h(x) = \frac{x}{\|x\|}$ for $\|x\| \geq R_0$, we have $d\xi_t = dt$ on $\{\|X_t\| \geq R_0\}$. Now we introduce an auxiliary process $Z$ which shall serve to control $\|X\|$.

Let $0 < \tilde{\eta} < \eta$. According to Skorokhod's lemma (see [16]), there is a unique pair of continuous adapted processes $(Z, L)$ such that $L$ is an increasing process (of finite variation) which increases only at times $t$ for which $Z_t = R_0$, and $Z \geq R_0$, which satisfies the equation

$$Z := R_0 \vee \|x\| + \sqrt{\varepsilon}M - \tilde{\eta}\xi + L.$$

We will prove that

(1.18) $$\|X_t\| \leq Z_t \quad \text{a.s. for all } t \geq 0.$$

For that purpose, choose $f \in C^2(\mathbb{R})$ such that

$$f(x) > 0 \quad \text{and} \quad f'(x) > 0 \quad \text{for all } x > 0,$$
$$f(x) = 0 \quad \text{for all } x \leq 0.$$

According to Itô's formula, for $t \geq 0$,

$$f(h(X_t) - Z_t) = f(h(x) - \|x\| \vee R_0) + \int_0^t f'(h(X_s) - Z_s) \, d(h(X) - Z)_s$$
$$+ \tfrac{1}{2} \int_0^t f''(h(X_s) - Z_s) \, d\langle h(X) - Z \rangle_s.$$

By definition of $h$ and $Z$, we have $h(X_t) \leq Z_t$ on $\{\|X_t\| \leq R_0\}$, so $\{h(X_t) > Z_t\} = \{\|X_t\| > Z_t\}$. Moreover, by definition, $h(X) - Z$ is a finite variation process. Hence, the expression

$$\int_0^t f'(\|X_s\| - Z_s)\left\{\frac{1}{\|X_s\|}\left\langle X_s, b\left(\frac{s}{T^\varepsilon}, X_s\right)\right\rangle + \frac{\varepsilon}{2}\triangle h(X_s) + \tilde{\eta}\right\} ds$$
$$- \int_0^t f'(\|X_s\| - Z_s) \, dL_s$$

is an upper bound of $f(h(X_t) - Z_t)$. Furthermore, $\triangle h(x) = \frac{d-1}{\|x\|}$ for $\|x\| \geq R_0$, which by (1.16) implies

$$\frac{1}{\|X_s\|}\left\langle X_s, b\left(\frac{s}{T^\varepsilon}, X_s\right)\right\rangle + \frac{\varepsilon}{2}\triangle h(X_s) + \tilde{\eta} < \frac{\varepsilon(d-1)}{2\|X_s\|} + \tilde{\eta} - \eta \quad \text{on } \{\|X_s\| > Z_s\}.$$

The latter expression is negative if $\varepsilon$ is small enough, so we can find some $\varepsilon_0 > 0$ such that $f(\|X_t\| - Z_t) \leq 0$ for $\varepsilon < \varepsilon_0$. This implies $\|X_t\| \leq Z_t$ a.s. by the definition of $f$, and (1.18) is established.

We therefore can bound the exit probability of $X$ by that of $Z$. If $Q$ denotes the law of the process $Z$, we see that, for any $\alpha > 0$,

(1.19) $$\mathbb{P}_x(\sigma_R \leq r(\varepsilon)) \leq Q(\sigma_R \leq r(\varepsilon)) \leq e^{\alpha r(\varepsilon)} \mathbb{E}_Q[e^{-\alpha \sigma_R}].$$

In order to find a bound on the right-hand side of (1.19), let us define $K := \sup_{\|x\| \leq R_0} \|\nabla h(x)\|^2$. Then we have $\xi_t \leq Kt$ for all $t \geq 0$. Note that

w.l.o.g. $h$ can be chosen so that $K \leq 2R_0$. Now observe that, by Itô's formula, for any $\varphi \in C^2(\mathbb{R})$,

$$d(\varphi(Z_t)e^{-(\alpha/K)\xi_t}) = \sqrt{\varepsilon}\varphi'(Z_t)e^{-(\alpha/K)\xi_t}\,dM_t + \varphi'(Z_t)e^{-(\alpha/K)\xi_t}\,dL_t$$
$$+ e^{-(\alpha/K)\xi_t}\left\{\frac{\varepsilon}{2}\varphi''(Z_t) - \tilde{\eta}\varphi'(Z_t) - \frac{\alpha}{K}\varphi(Z_t)\right\}d\xi_t.$$

Now let $R \geq R_0$. If we choose $\varphi$ such that

$$\frac{\varepsilon}{2}\varphi''(y) - \tilde{\eta}\varphi'(y) - \frac{\alpha}{K}\varphi(y) = 0 \qquad \text{for } y \in [R_0, R],$$
$$\varphi'(R_0) = 0, \qquad \varphi(R) = 1,$$

then $\varphi(Z_t)e^{-(\alpha/K)\xi_t}$ is a local martingale which is bounded up to time $\sigma_R$. Hence, we are allowed to apply the stopping theorem to obtain

$$(1.20) \qquad \varphi(R_0 \vee \|x\|) = \mathbb{E}_Q[\varphi(Z_{\sigma_R})e^{-(\alpha/K)\xi_{\sigma_R}}] = \mathbb{E}_Q[e^{-(\alpha/K)\xi_{\sigma_R}}].$$

Hence, $\xi_{\sigma_R} \leq K\sigma_R$, which implies $\mathbb{E}_Q[e^{-(\alpha/K)\xi_{\sigma_R}}] \geq \mathbb{E}_Q[e^{-\alpha\sigma_R}]$, and we deduce from (1.19) that

$$(1.21) \qquad \mathbb{P}_x(\sigma_R \leq r(\varepsilon)) \leq e^{\alpha r(\varepsilon)}\mathbb{E}_Q[e^{-(\alpha/K)\xi_{\sigma_R}}] \leq e^{\alpha r(\varepsilon)}\varphi(R_0 \vee \|x\|).$$

Solving the differential equation for $\varphi$ yields

$$\varphi(x) = \frac{-\lambda^- e^{\lambda^+(x-R_0)} + \lambda^+ e^{\lambda^-(x-R_0)}}{-\lambda^- e^{\lambda^+(R-R_0)} + \lambda^+ e^{\lambda^-(R-R_0)}},$$

with $\lambda^\pm = \frac{\tilde{\eta} \pm \sqrt{\tilde{\eta}^2 + 2(\alpha/K)\varepsilon}}{\varepsilon}$. Hence,

$$\varphi(x) \leq \frac{(\lambda^+ - \lambda^-)e^{\lambda^+(x-R_0)}}{(-\lambda^-)e^{\lambda^+(R-R_0)}}.$$

Taking $\alpha = r(\varepsilon)^{-1}$ in (1.21), we obtain

$$\mathbb{P}_x(\sigma_R \leq r(\varepsilon)) \leq \exp(1)\varphi(R_0 \vee \|x\|) \leq \frac{\lambda^+ - \lambda^-}{-\lambda^-}\exp\{1 + \lambda^+(R_0 \vee \|x\| - R)\}.$$

It is obvious that $\exp\{\lambda^+(R_0 \vee \|x\| - R)\} \leq \exp\{-\frac{\tilde{\eta}R}{\varepsilon}\}$ for $R \geq 2(\|x\| \vee R_0)$, so it remains to comment on the prefactor. We have

$$\frac{\lambda^+ - \lambda^-}{-\lambda^-} = \frac{2\sqrt{\tilde{\eta}^2 + 2(\alpha/K)\varepsilon}}{\sqrt{\tilde{\eta}^2 + 2(\alpha/K)\varepsilon} - \tilde{\eta}} \leq \frac{4(\tilde{\eta}^2 + 2\varepsilon/(Kr(\varepsilon)))}{2\varepsilon/(Kr(\varepsilon))}.$$

Since $\lim_{\varepsilon \to 0}\frac{\varepsilon}{r(\varepsilon)} = 0$, the latter expression behaves like $2\tilde{\eta}^2 K\frac{r(\varepsilon)}{\varepsilon}$ as $\varepsilon \to 0$. Putting these estimates together yields the claimed asymptotic bound with $\tilde{\eta}$ instead of $\eta$. Letting $\tilde{\eta} \to \eta$ completes the proof. $\square$

REMARK 1.5. The proof of Proposition 1.4 shows a lot more. The crucial inequality (1.21) contains a bound which is independent of $X^\varepsilon$, since $\varphi$ is defined by means of $h$, $\varepsilon$, $\tilde\eta$ and $R_0$ only. Thus, we have shown that the bound (1.17) holds for *all* diffusions satisfying the growth condition (1.16), that is, $\varepsilon_1$ and $R_1$ are independent of $X^\varepsilon$. In particular, (1.17) holds uniformly w.r.t. $\mu$.

1.3.2. *Properties of the quasi-potential.* Taking large period limits in the subsequently derived large deviations results for our diffusions with slow periodic variation will require to freeze the time parameter in the drift term. The corresponding rate functions are given a separate treatment in this subsection. We shall briefly discuss their regularity properties. This will be of central importance for the estimation of exit rates in Section 2. For $s \geq 0, T > 0$, we consider

$$(1.22) \quad I_{0T}^s(\varphi) = \begin{cases} \frac{1}{2} \int_0^T \|\dot\varphi_t - b(s, \varphi_t)\|^2 \, dt, \\ \qquad \text{if } \varphi \text{ is absolutely continuous,} \\ +\infty, \qquad \text{otherwise.} \end{cases}$$

As in the first section, we need associated *cost* functions. For $s \geq 0$, $x, y \in \mathbb{R}^d$, they are given by

$$(1.23) \quad V^s(x, y, t) = \inf\{I_{0t}^s(\varphi) : \varphi \in C_{0t}, \varphi_0 = x, \varphi_t = y\}.$$

$V^s(x, y, t)$ is the cost of forcing the *frozen* system

$$dY_t^\varepsilon = b(s, Y_t^\varepsilon) \, dt + \sqrt{\varepsilon} \, dW_t, \qquad t \geq 0,$$

to be at the point $y$ at time $t$ when starting at $x$. The corresponding *quasi-potential*

$$(1.24) \quad V^s(x, y) = \inf_{t > 0} V^s(x, y, t)$$

describes the cost for the frozen system to go from $x$ to $y$ eventually. Let us note that since the drift $b$ is locally Lipschitz in the time variable, the family of action functionals $I_{0T}^s$ is continuous w.r.t. the parameter $s$, and the corresponding cost functions and pseudo-potentials inherit this continuity property. Let us recall some further useful properties of the quasi-potentials and their underlying cost and rate functions. The following properties are immediate.

LEMMA 1.6. *For any $x, y, z \in \mathbb{R}^d$ and $s, t, u \geq 0$, we have the following:*

(a) $V^s(x, y, t + u) \leq V^s(x, z, t) + V^s(z, y, u)$,
(b) $(s, y) \mapsto V^s(x, y, t)$ *is continuous on* $\mathbb{R}_+ \times \mathbb{R}^d$,
(c) $\inf_{\|y\| \geq R} V^s(x, y, t) \xrightarrow[R \to \infty]{} \infty$ *uniformly w.r.t.* $s \geq 0$.

The following lemma establishes the local Lipschitz continuity of the quasi-potential w.r.t. the state variables, uniformly w.r.t. the parameter $s$.

LEMMA 1.7. *For any compact subset $K$ of $\mathbb{R}^d$, there exists $\Gamma_K \geq 0$ such that*

$$\sup_{s \geq 0} V^s(x, y) \leq \Gamma_K \operatorname{dist}(x, y)$$

*for all $x, y \in K$.*

PROOF. Let $x$ and $y$ belong to $K$. There exists some radius $R > 0$ such that $K \subset B_R(0)$. Set $T = \operatorname{dist}(x, y)$. We construct a path $\varphi \in C_{0T}$ by setting $\varphi_t = x + \frac{y-x}{\operatorname{dist}(x,y)} t$ for $t \in [0, T]$. Since $b(s, \cdot)$ is locally Lipschitz uniformly w.r.t. $s \geq 0$, we obtain an upper bound for the energy of $\varphi$:

$$\begin{aligned}
I_{0T}^s(\varphi) &\leq \frac{1}{2} \sup_{u \geq 0} \int_0^T \|\dot{\varphi}_t - b(u, \varphi_t)\|^2 \, dt \\
&\leq \frac{1}{2} \int_0^T \left( \frac{\|y - x\|}{\operatorname{dist}(x, y)} + \sup_{0 \leq u \leq 1} \|b(u, \varphi_t)\| \right)^2 dt \\
&\leq \frac{1}{2} \int_0^T (1 + \kappa_R(0) + \|b(0, \varphi_t)\|)^2 \, dt \\
&\leq \frac{1}{2} \int_0^T (1 + \kappa_R(0) + K_R(0)\|\varphi_t\| + \|b(0, 0)\|)^2 \, dt \\
&\leq \frac{T}{2} (1 + \kappa_R(0) + R K_R(0) + \|b(0, 0)\|)^2.
\end{aligned}$$

For $\Gamma_K := \frac{1}{2}(1 + \kappa_R(0) + R K_R(0) + \|b(0, 0)\|)^2$ and by the definition of $T$, we obtain

$$\sup_{s \geq 0} V^s(x, y) \leq \sup_{s \geq 0} I_{0T}^s(\varphi) \leq \Gamma_K \operatorname{dist}(x, y). \qquad \square$$

1.3.3. *Large deviations.* We shall now specialize the general large deviations results of the previous subsection to the family $X^\varepsilon$, $\varepsilon > 0$, of solutions of (1.12). At the same time they will be strengthened, to obtain uniformity w.r.t. some of the system's parameters: the scale parameter $\mu$, the starting time and the initial condition.

It is an immediate consequence of Proposition 1.3 that the solution of (1.12) satisfies a large deviations principle with rate function $I_{0T}^0$, that is, the rate function is the same as that of a homogeneous diffusion governed by the frozen drift $b(0, \cdot)$. In order to see this, one only has to mention that $\lim_{\varepsilon \to 0} b(\frac{t}{T\varepsilon}, x) = b(0, x)$ locally uniformly w.r.t. $x$ due to the Lipschitz assumptions on $b$.

But this result is not strong enough. We also need some uniformity w.r.t. the starting times of the diffusions we consider. Our large deviations statements derived so far rely on comparison arguments which yield exponential equivalence with time homogeneous diffusions for which an LDP is well known from the classical theory of Freidlin and Wentzell. In order to achieve uniform large deviations estimates, we have to refine this technique, to derive a large deviations principle for our family of diffusions (1.12), which is uniform with respect to both the starting time and the scale parameter. This will be our main tool for estimating the asymptotics of exit time laws in the subsequent section.

The diffusion (1.12) is a time inhomogeneous Markov process. The solution starting at time $r \geq 0$ with initial condition $x \in \mathbb{R}^d$ has the same law as the solution $X^{r,x}$ of the SDE

$$(1.25) \quad dX_t^{r,x} = b\left(\frac{r+t}{T^\varepsilon}, X_t^{r,x}\right) dt + \sqrt{\varepsilon}\, dW_t, \qquad t \geq 0, X_0^{r,x} = x \in \mathbb{R}^d.$$

We denote its law by $\mathbb{P}_{x,r}(\cdot)$, assume from now on that $T^\varepsilon = \exp \frac{\mu}{\varepsilon}$ for some $\mu > 0$, and fix $T \geq 0$.

PROPOSITION 1.8. *Let $K \subset \mathbb{R}^d$ be a compact set and $\mathcal{V} \subset (0,\infty)$. For $\mu \in \mathcal{V}$, $r \in [0,1]$ and $\beta \geq 0$, let $S^{r,\beta}(\varepsilon,\mu)$ be a neighborhood of $rT^\varepsilon$ such that*

$$\limsup_{\varepsilon \to 0} \sup_{\mu \in \mathcal{V}, r \in [0,1]} \frac{\mathrm{diam}(S^{r,\beta}(\varepsilon,\mu))}{T^\varepsilon} \leq \beta.$$

*Then for any closed $F \subset C_{0T}$, there exists $\delta = \delta(F)$ such that*

$$\limsup_{\varepsilon \to 0} \varepsilon \log \sup_{y \in K, \mu \in \mathcal{V}, u \in S^{r,\beta}(\varepsilon,\mu)} \mathbb{P}_{y,u}(X^\varepsilon \in F) \leq -\inf_{y \in K} \inf_{\varphi \in F^{\gamma_0}, \varphi_0 = y} I_{0T}^r(\varphi),$$

*where $\gamma_0 = \gamma_0(F) = \beta\delta(F)$ and $F^{\gamma_0}$ is the closed $\gamma_0$-neighborhood of $F$. For any open $G \subset C_{0T}$, there exists $\delta = \delta(G)$ and $\beta_0 = \beta_0(G)$ such that, if $\beta \leq \beta_0$,*

$$\liminf_{\varepsilon \to 0} \varepsilon \log \inf_{y \in K, \mu \in \mathcal{V}, u \in S^{r,\beta}(\varepsilon,\mu)} \mathbb{P}_{y,u}(X^\varepsilon \in G) \geq -\sup_{y \in K} \inf_{\varphi \in G^{\gamma_0}, \varphi_0 = y} I_{0T}^r(\varphi),$$

*where $\gamma_0 = \gamma_0(G) = \beta\delta(G)$ and $G^{\gamma_0}$ is the complement of $(G^c)^{\gamma_0}$.*

*These bounds hold uniformly w.r.t. $r$.*

REMARK 1.9. The upper bound means that, for any $\vartheta > 0$, we can find $\varepsilon_0 > 0$ s.t. for $\varepsilon \leq \varepsilon_0$, we have

$$\varepsilon \log \sup_{y \in K, \mu \in \mathcal{V}, u \in S^{r,\beta}(\varepsilon,\mu)} \mathbb{P}_{y,u}(X^\varepsilon \in F) \leq -\inf_{y \in K} \inf_{\varphi \in F^{\gamma_0}, \varphi_0 = y} I_{0T}^r(\varphi) + \vartheta.$$

The uniformity in the statement means that $\varepsilon_0$ can be chosen independently of $r$. A similar statement holds for the lower bound.

Observe that the expression for the blowup-factor $\gamma_0(F)$ depends on the set $F$ only through $\delta(F)$ which is independent of $\beta$, and that $\gamma_0(F) \to 0$ as $\beta \to 0$ for all $F$. In particular, if $\beta$ is equal to zero, we recover the classical bound of the uniform LDP.

PROOF OF PROPOSITION 1.8. For $y \in \mathbb{R}^d$ and $r \in [0,1]$, let $Y^{r,y}$ be the solution of the homogeneous SDE

$$dY_t^{r,y} = b(r, Y_t^{r,y})\, dt + \sqrt{\varepsilon}\, dW_t, \qquad t \geq 0, Y_0^{r,y} = y.$$

Let $\mathcal{W} \subset [0,1]$ and $r_0 \in \mathcal{W}$. For $r \in \mathcal{W}$, $u \in S^{r,\beta}(\varepsilon,\mu), \mu \in \mathcal{V}$ and $R > 0$, let $\tau_R^{u,y} := \inf\{t \geq 0 : X_t^{u,y} \notin B_R(0)\}$, and let $\tilde{\tau}_R^{r_0,y}$ be defined similarly with $X^{u,y}$ replaced by $Y^{r_0,y}$, and $\sigma_R^{u,y,r_0} := \tau_R^{u,y} \wedge \tilde{\tau}_R^{r_0,y}$.

As a consequence of Gronwall's lemma, we see just as in the proof of Proposition 1.3 that, for $r, r_0 \in [0,1], u \in S^{r,\beta}(\varepsilon,\mu)$,

$$\rho_{0T}(X^{u,y}, Y^{r_0,y}) \leq e^{K_R(0)T} \int_0^T \left\| b\left(\frac{u+t}{T^\varepsilon}, X_t^{u,y}\right) - b(r_0, X_t^{u,y}) \right\| dt$$

on $\{\sigma_R^{u,y,r_0} > T\}$. This implies

$$\rho_{0T}(X^{u,y}, Y^{r_0,y}) \leq \kappa_R(0) T e^{K_R(0)T} \left( \frac{\operatorname{diam}(S^{r,\beta}(\varepsilon,\mu)) + T}{T^\varepsilon} + |r - r_0| \right)$$

on $\{\sigma_R^{u,y,r_0} > T\}$. Due to our assumption, the last expression is bounded by

(1.26) $\qquad \beta_1 = \beta_1(\mathcal{W}) = \beta_0(\mathcal{W}) M(R) \qquad$ as $\varepsilon \to 0$,

where

(1.27) $\quad \beta_0(\mathcal{W}) := \beta + \sup_{r \in \mathcal{W}} |r - r_0| \quad$ and $\quad M(R) := T \kappa_R(0) e^{K_R(0)T}$.

*Upper bound.* Fix some closed set $F \subset C_{0T}$. For all $\gamma > 0$, we have

$$\mathbb{P}(X^{u,y} \in F) \leq \mathbb{P}(Y^{r_0,y} \in F^\gamma) + \mathbb{P}(\rho_{0T}(X^{u,y}, Y^{r_0,y}) > \gamma).$$

This yields

$$\limsup_{\varepsilon \to 0} \varepsilon \log \sup_{y \in K, \mu \in \mathcal{V}, r \in \mathcal{W}, u \in S^{r,\beta}(\varepsilon,\mu)} \mathbb{P}_{y,u}(X^\varepsilon \in F)$$

(1.28) $\quad \leq \limsup_{\varepsilon \to 0} \varepsilon \log \max \Big\{ \sup_{y \in K} \mathbb{P}(Y^{r_0,y} \in F^\gamma),$

$$\sup_{y \in K, \mu \in \mathcal{V}, r \in \mathcal{W}, u \in S^{r,\beta}(\varepsilon,\mu)} \mathbb{P}(\rho_{0T}(X^{u,y}, Y^{r_0,y}) > \gamma) \Big\}.$$

Now we wish to find some $\gamma$ such that the asymptotics of the maximum is determined by the left term $\sup_{y\in K}\mathbb{P}(Y^{r_0,y}\in F^\gamma)$. In that case the uniform LDP for $Y$ will give us the bound

$$(1.29) \quad \limsup_{\varepsilon\to 0}\varepsilon\log\sup_{y\in K,\mu\in\mathcal{V},r\in\mathcal{W},u\in S^{r,\beta}(\varepsilon,\mu)}\mathbb{P}_{y,u}(X^\varepsilon\in F)$$
$$\leq -\inf_{y\in K}\inf_{\varphi\in F^\gamma,\varphi_0=y}I_{0T}^{r_0}(\varphi).$$

Unfortunately, such a $\gamma$ will depend on $F$. In order to see that it exists and can be chosen as claimed in the statement, we define

$$\Theta(R,\varepsilon):=\sup_{r\in[0,1],y\in K,\mu\in\mathcal{V},u\in S^{r,\beta}(\varepsilon,\mu)}\mathbb{P}(\tau_R^{u,y}\leq T)+\sup_{r\in[0,1],y\in K}\mathbb{P}(\tilde{\tau}_R^{r,y}\leq T).$$

By Proposition 1.4 and Remark 1.5, we have $\limsup_{\varepsilon\to 0}\varepsilon\log\Theta(R,\varepsilon)\leq -\eta R$ for all $R\geq R_1$. Hence, we may find $R\geq R_1$ such that

$$\limsup_{\varepsilon\to 0}\varepsilon\log\Theta(R,\varepsilon)\leq -\sup_{r\in[0,1]}\inf_{y\in K}\inf_{\varphi\in F,\varphi_0=y}I_{0T}^r(\varphi).$$

Fix this $R$, let $\delta(F)=M(R)$, and note that $\delta(F)$ is independent of $\beta$ and $\mathcal{W}$. By (1.26), for any $\gamma>\beta_1(\mathcal{W})=\beta_0(\mathcal{W})\delta(F)$, we can find $\varepsilon_0>0$ such that, for $\varepsilon\leq\varepsilon_0$,

$$(1.30) \quad \sup_{r\in\mathcal{W},y\in K,\mu\in\mathcal{V},u\in S^{r,\beta}(\varepsilon,\mu)}\mathbb{P}(\rho_{0T}(X^{u,y},Y^{r_0,y})>\gamma)\leq\Theta(R,\varepsilon).$$

Now we exploit the definition of $\delta(F)$, which ensures that the maximum in (1.28) is given by the left term. Indeed, for $\gamma>\beta_1(\mathcal{W})$, we have

$$\limsup_{\varepsilon\to 0}\varepsilon\log\sup_{r\in\mathcal{W},y\in K,\mu\in\mathcal{V},u\in S^{r,\beta}(\varepsilon,\mu)}\mathbb{P}(\rho_{0T}(X^{u,y},Y^{r_0,y})>\gamma)$$
$$\leq\limsup_{\varepsilon\to 0}\varepsilon\log\Theta(R,\varepsilon)\leq -\sup_{r\in[0,1]}\inf_{y\in K}\inf_{\varphi\in F,\varphi_0=y}I_{0T}^r(\varphi)$$
$$\leq -\inf_{y\in K}\inf_{\varphi\in F^\gamma,\varphi_0=y}I_{0T}^{r_0}(\varphi),$$

which implies (1.29). The particular choice $\mathcal{W}=\{r_0\}$ yields this bound for all $\gamma>\gamma_0(F)=\beta\delta(F)$ given in the statement and proves the claimed bound. By taking the limit $\gamma\to\gamma_0(F)$, we get the asserted upper bound since $I$ is a good rate function (see [4], Lemma 4.1.6).

It remains to prove the uniformity w.r.t. $r$. For that purpose, fix $\vartheta>0$, and for $r_0\in[0,1]$, let $\mathcal{W}_{r_0}$ be a neighborhood of $r_0$. Choose $\phi^*\in F^{\gamma_0}$ [here $\gamma_0=\gamma_0(F)=\beta\delta(F)$] starting at some point $y_0\in K$ with $\inf_{y\in K}\inf_{\varphi\in F^{\gamma_0},\varphi_0=y}I_{0T}^{r_0}(\varphi)\geq I_{0T}^{r_0}(\phi^*)-\vartheta/8$. By Lemma 4.1.6 in [4] and the continuity of $r\mapsto I_{0T}^r(\phi^*)$, we can assume $\mathcal{W}_{r_0}$ to be small enough such that, for $r \in \mathcal{W}_{r_0}$, denoting $\gamma^* = \beta_0(\mathcal{W}_{r_0})\delta(F)$,

$$\inf_{y \in K} \inf_{\varphi \in F^{\gamma^*}, \varphi_0 = y} I_{0T}^{r_0}(\varphi) \geq \inf_{y \in K} \inf_{\varphi \in F^{\gamma_0}, \varphi_0 = y} I_{0T}^{r_0}(\varphi) - \vartheta/8$$
$$\geq I_{0T}^{r_0}(\phi^*) - \vartheta/4$$
$$\geq I_{0T}^{r}(\phi^*) - \vartheta/2$$
$$\geq \inf_{y \in K} \inf_{\varphi \in F^{\gamma_0}, \varphi_0 = y} I_{0T}^{r}(\varphi) - \vartheta/2.$$

Due to compactness, we can choose finitely many points $r_1, ..., r_N$ such that their corresponding neighborhoods cover $[0, 1]$. Denote $\gamma_n^* := \beta_0(\mathcal{W}_{r_n})\delta(F)$. For each $1 \leq n \leq N$, there exists some $\varepsilon_n > 0$ such that, for $\varepsilon \leq \varepsilon_n$ and $r \in \mathcal{W}_{r_n}$,

$$\varepsilon \log \sup_{y \in K, \mu \in \mathcal{V}, u \in S^r(\varepsilon, \mu)} \mathbb{P}_{y,u}(X^\varepsilon \in F) \leq -\inf_{y \in K} \inf_{\varphi \in F^{\gamma_n^*}, \varphi_0 = y} I_{0T}^{r_n}(\varphi) + \frac{\vartheta}{2}$$
$$\leq -\inf_{y \in K} \inf_{\varphi \in F^{\gamma_0}, \varphi_0 = y} I_{0T}^{r}(\varphi) + \vartheta.$$

Hence, for $\varepsilon \leq \min_{1 \leq n \leq N} \varepsilon_n$, the preceding inequality holds for all $r \in [0, 1]$.

*Lower bound.* Let $G \subset C_{0T}$ be an open set. Consider the increasing function

$$f(l) := \frac{1}{\eta} \sup_{y \in K} \inf_{\phi \in G^l : \phi_0 = y} I_{0T}^{r_0}(\phi),$$

let $l_0 = \inf\{l \geq 0 : f(l) = +\infty\}$, and recall that $\eta$ is the constant introduced in the growth condition (1.16) for the drift. Assume first that $l_0 < \infty$ (this is guaranteed if $G$ is bounded), and set

$$R := f\left((l_0 - \beta_0(\mathcal{W})) \vee \frac{l_0}{2}\right) \vee R_1 \quad \text{and} \quad \gamma := \beta_0(\mathcal{W}) M(R),$$

where $R_1$ is given by Proposition 1.4. Then

$$\mathbb{P}(Y^{r_0, y} \in G^\gamma) \leq \mathbb{P}(X^{u,y} \in G) + \mathbb{P}(\rho_{0T}(Y^{r_0, y}, X^{u,y}) > \gamma).$$

By the uniform LDP for $Y^{r_0, y}$ and (1.30), we conclude that

$$-\eta f(\gamma) = -\sup_{y \in K} \inf_{\varphi \in G^\gamma, \varphi_0 = y} I_{0T}^{r_0}(\varphi) \leq \liminf_{\varepsilon \to 0} \varepsilon \log \inf_{y \in K} \mathbb{P}(Y^{r_0, y} \in G^\gamma)$$
$$\leq \max\left\{\liminf_{\varepsilon \to 0} \varepsilon \log \inf_{r \in \mathcal{W}, y \in K, \mu \in \mathcal{V}, u \in S^r(\varepsilon, \mu)} \mathbb{P}(X^{u,y} \in G),\right.$$
$$\left. \limsup_{\varepsilon \to 0} \varepsilon \log \sup_{r \in \mathcal{W}, y \in K, \mu \in \mathcal{V}, u \in S^r(\varepsilon, \mu)} \mathbb{P}(\rho_{0T}(Y^{r_0, y}, X^{u,y}) > \gamma)\right\}$$

$$\leq \max\left\{\liminf_{\varepsilon \to 0} \varepsilon \log \inf_{r \in \mathcal{W}, y \in K, \mu \in \mathcal{V}, u \in S^r(\varepsilon, \mu)} \mathbb{P}(X^{u,y} \in G),\right.$$
$$\left.\limsup_{\varepsilon \to 0} \varepsilon \log \Theta(R, \varepsilon)\right\}.$$

Since $f$ is increasing and $R \geq R_1$, we obtain by Proposition 1.4

$$-\eta f(\gamma + \beta_0(\mathcal{W}))$$
$$\leq -\eta f(\gamma)$$
$$\leq \max\left\{\liminf_{\varepsilon \to 0} \varepsilon \log \inf_{r \in \mathcal{W}, y \in K, \mu \in \mathcal{V}, u \in S^r(\varepsilon, \mu)} \mathbb{P}(X^{u,y} \in G), -\eta R\right\}.$$

Now we have to compare $f(\gamma)$ and $R$ in order to see when the maximum is given by the left term.

If $f(\gamma) > R$, then $\gamma > (l_0 - \beta_0(\mathcal{W})) \vee \frac{l_0}{2} \geq l_0 - \beta_0(\mathcal{W})$ by monotonicity of $f$, hence, $f(\gamma + \beta_0(\mathcal{W})) = +\infty$ by definition of $l_0$. Otherwise, we have $f(\gamma) \leq R$, which means that the left term dominates the maximum.

In both cases we get

$$-\eta f(\gamma + \beta_0(\mathcal{W})) \leq \liminf_{\varepsilon \to 0} \varepsilon \log \inf_{r \in \mathcal{W}, y \in K, \mu \in \mathcal{V}, u \in S^r(\varepsilon, \mu)} \mathbb{P}(X^{u,y} \in G).$$

Now consider the unbounded case $l_0 = +\infty$. Recall $M$ defined by (1.27), and let $\beta_0(G) := \sup_{l \geq 0} \frac{l}{M(f(l))}$, the existence of which was claimed in the statement. If $\beta_0(\mathcal{W}) < \beta_0(G)$, we can choose $l_1$ such that $\frac{l_1}{M(f(l_1))} \geq \beta_0(\mathcal{W})$ and set $\gamma := \beta_0(\mathcal{W}) M(f(l_1))$. Using the same arguments as in the bounded case, we deduce that

$$-\eta f(\gamma) \leq \max\left\{\liminf_{\varepsilon \to 0} \varepsilon \log \inf_{r \in \mathcal{W}, y \in K, \mu \in \mathcal{V}, u \in S^r(\varepsilon, \mu)} \mathbb{P}(X^{u,y} \in G), -\eta f(l_1)\right\}.$$

Since $f$ is increasing and $l_1 \geq \gamma$, we obtain

$$-\eta f(\gamma) \leq \liminf_{\varepsilon \to 0} \varepsilon \log \inf_{r \in \mathcal{W}, y \in K, \mu \in \mathcal{V}, u \in S^r(\varepsilon, \mu)} \mathbb{P}(X^{u,y} \in G).$$

In both the bounded and the unbounded case we have found $\gamma = \beta_0(\mathcal{W})\delta(G)$ such that the desired bound holds: we have $\delta(G) = M(R) + 1$ in the bounded case and $\delta(G) = M(f(l_1))$ in the unbounded case. Furthermore, the choise $\mathcal{W} = \{r_0\}$ corresponds to $\beta_0(\mathcal{W}) = \beta$ and yields $\gamma_0(G) = \beta\delta(G)$, in complete analogy to the situation of the upper bound. The uniformity is also proved in exactly the same way as already shown for the upper bound. $\square$

**2. Exit and entrance times of domains of attraction.** We continue to study asymptotic properties of diffusions with weakly periodic drifts given by the SDE

$$(2.1) \qquad dX^\varepsilon_t = b\left(\frac{t}{T^\varepsilon}, X^\varepsilon_t\right) dt + \sqrt{\varepsilon}\, dW_t, \qquad t \geq 0, X^\varepsilon_0 = x_0 \in \mathbb{R}^d.$$

In this section we shall work out the effects of weak periodicity of the drift on the asymptotic behavior of the exit times of its domains of attraction. This will be done under simple assumptions on the geometry associated to it. So we will have to specify some assumptions on the attraction and conservation properties of $b$. Essentially, we shall assume that $\mathbb{R}^d$ is split into two domains of attraction, separated by a simple geometric boundary which is invariant in time. Apart from that, we shall assume that the drift is pointing inward sufficiently strongly so that the diffusions will not be able to leave compact sets in the small noise limit. Let us make these assumptions more precise. We recall that, according to the Kramers–Eyring law (see, e.g., [9]), the mean time a homogeneous diffusion of noise intensity $\varepsilon$ needs to leave a potential well of depth $\frac{v}{2}$ is of the order $\exp \frac{v}{\varepsilon}$. Nature therefore imposes the time scales $T^\varepsilon$ with which we have to work. For simplicity, we measure these scales in energy units: with $\mu > 0$, we associate the time scale $T^\varepsilon = \exp \frac{\mu}{\varepsilon}$. We assume as before that $b$ satisfies the local Lipschitz conditions (1.14) and (1.15), and that the growth of the inward drift is sufficiently strong near infinity which is expressed by (1.16).

The additional conditions concerning the geometry of $b$ are specified in the following.

ASSUMPTION 2.1. The two-dimensional ordinary differential equation

$$\dot\varphi_s(t) = b(s, \varphi_s(t)), \qquad t \geq 0, \tag{2.2}$$

admits two stable equilibria $x_-$ and $x_+$ in $\mathbb{R}^d$ which do not depend on $s \geq 0$. Moreover, the domains of attraction defined by

$$A_\pm(s) = \left\{ y \in \mathbb{R}^d : \dot\varphi_s(t) = b(s, \varphi_s(t)) \right.$$
$$\left. \text{and } \varphi_s(0) = y \text{ imply } \lim_{t \to \infty} \varphi_s(t) = x_\pm \right\} \tag{2.3}$$

are also independent of $s \geq 0$ and denoted by $A_\pm$. They are supposed to satisfy $\overline{A_- \cup A_+} = \mathbb{R}^d$, and $\partial A_- = \partial A_+$. We denote by $\chi$ the common boundary (see Figure 1).

The asymptotic results concerning the exit and entrance time should remain true in a more general setting, where the stable equilibrium points and the domains of attraction depend on $s \geq 0$. We stick to Assumption 2.1 for reasons of technical and notational simplicity.

The main subject of investigation in this section is given by the exit times of the domains of attraction $A_\pm$, provided that the weakly time inhomogeneous diffusion starts near the equilibrium points $x_\pm$. By obvious symmetry reasons, we may restrict our attention to the case of an exit from $A_-$. As we shall show, this exit time depends on the quasi-potential, that is, on the cost function taken on the set of all functions starting in the neighborhood of $x_-$ and exiting the domain of attraction through $\chi$. For this reason, we introduce the one-periodic *energy function*

$$(2.4) \qquad e(s) := \inf_{y \in \chi} V^s(x_-, y) < \infty \qquad \text{for } s \geq 0,$$

which is continuous on $\mathbb{R}_+$. In the gradient case $b(t,x) = -\nabla_x U(t,x)$, this function coincides with twice the depth of the potential barrier to be overcome in order to exit from $A_-$, that is, the energy the diffusion needs to leave $A_-$. Therefore, scales $\mu$—corresponding to the Kramers–Eyring times $T^\varepsilon = \exp(\frac{\mu}{\varepsilon})$ according to the chosen parametrization—at which we expect transitions between the domains of attraction must be comprised between

$$\mu_* := \inf_{t \geq 0} e(t) \quad \text{and} \quad \mu^* := \sup_{t \geq 0} e(t).$$

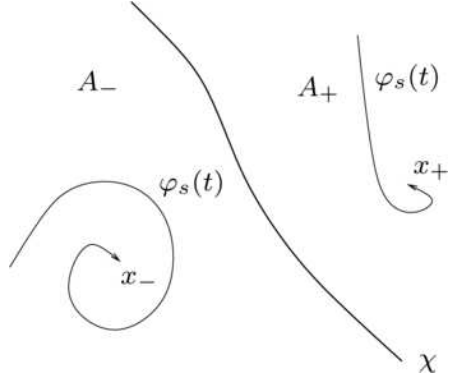

Fig. 1. *Domains of attraction.*

These two constants are finite and are reached at least once per period since $e(t)$ is continuous and periodic. Now fix a time scale parameter $\mu$. This parameter serves as a threshold for the energy, and we expect to observe an exit from $A_-$ at the first time $t$ at which $e(t)$ falls below $\mu$. For $\mu \in ]\mu_*, \mu^*[$, we therefore define (see Figure 2).

$$(2.5) \qquad a_\mu = \inf\{t \geq 0 : e(t) \leq \mu\}, \qquad \alpha_\mu = \inf\{t \geq 0 : e(t) < \mu\}.$$

The subtle difference between $a_\mu$ and $\alpha_\mu$ may be important, but we shall rule it out for our considerations by making the following assumption.

ASSUMPTION 2.2. The energy function $e(t)$ is strictly monotonous between its (discrete) extremes, and every local extremum is global.

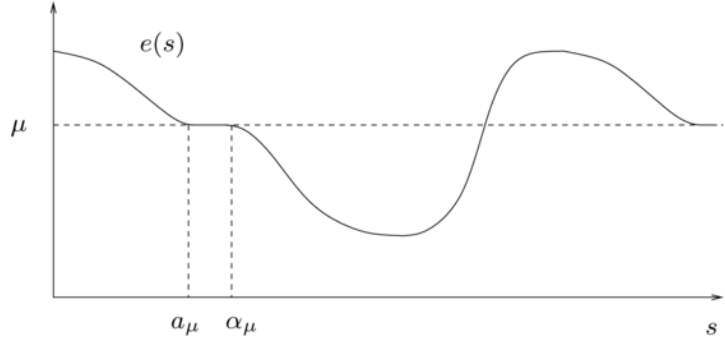

Fig. 2. *Definition of $a_\mu$ and $\alpha_\mu$.*

Under this assumption, we have $a_\mu = \alpha_\mu$. We are now in a position to state the main result of this section. Let $\varrho > 0$ be small enough such that the Euclidean ball $B_\varrho(x_+) \subset A_+$, and let us define the first entrance time into this ball by

$$\tau_\varrho = \inf\{t \geq 0 : X_t^\varepsilon \in B_\varrho(x_+)\}. \tag{2.6}$$

This stopping time depends of course on $\varepsilon$, but, for notational simplicity, we suppress this dependence.

THEOREM 2.3. *Let $\mu < e(0)$. There exist $\eta > 0$ and $h_0 > 0$ such that, for $h \leq h_0$,*

$$\lim_{\varepsilon \to 0} \varepsilon \log \sup_{y \in B_\eta(x_-)} \mathbb{P}_y(\tau_\varrho \notin [(a_\mu - h)T^\varepsilon, (\alpha_\mu + h)T^\varepsilon]) = \mu - e(a_\mu - h).$$

*Moreover, under Assumption 2.2, this convergence is uniform w.r.t. $\mu$ on compact subsets of $]\mu_*, e(0)[$.*

Note that Assumption 2.2 implies the continuity of $\mu \mapsto \mu - e(a_\mu - h)$, and that $\mu - e(a_\mu - h) < 0$ if $h$ is small enough. The statement of the theorem may be paraphrased in the following way. It specifies time windows in which transitions between the domains of attraction will be observed with very high probability. In particular, if $e(t)$ is strictly monotonous between its extremes, we prove that the entrance time into a neighborhood of $x_+$ will be located near $a_\mu T^\varepsilon$ in the small noise limit. The assumption $\mu < e(0)$ is only a technical assumption in order to avoid instantaneous jumping of the diffusion to the other valley. It can always be achieved by simply starting the diffusion a little later. We could even assume that $e(0) = \mu^*$, which then would yield uniform convergence on compact subsets of $]\mu_*, \mu^*[$.

The rest of this section is devoted to the proof of this main result and is subdivided into separate subsections in which lower and upper bounds are established.

2.1. *Lower bound for the exponential exit rate*: *diffusion exit.* We have to establish upper and lower bounds on the transition time $\tau_\varrho$ which both should be exceeded with an exponentially small probability that has to be determined exactly. It turns out that the probability of exceeding the upper bound $(\alpha_\mu - h)T^\varepsilon$ vanishes asymptotically to all exponential orders, so the exact large deviations rate is determined only by the probability $\mathbb{P}_x(\tau_\varrho \leq (a_\mu - h)T^\varepsilon)$ of exceeding the lower bound.

For a lower bound of the latter probability, as well as for an upper bound on $\mathbb{P}_x(\tau_\varrho \geq (\alpha_\mu + h)T^\varepsilon)$, one has to prove large deviations type upper bounds of the asymptotic distribution $\mathbb{P}_x(\tau_\varrho \geq s(\varepsilon))$ for suitably chosen $s(\varepsilon)$. This can be expressed in terms of the problem of *diffusion exit* from a carefully chosen bounded domain.

Recall that $\tau_\varrho$ is the *first entrance time* of a small neighborhood $B_\varrho(x_+)$ of the equilibrium point $x_+$. Consider for $R, \varrho > 0$ the bounded domain

$$D = D(R, \varrho) := \overline{B_R(0)} \setminus B_\varrho(x_+),$$

and let

$$\tau_D := \inf\{t \geq 0 : X_t \notin D\}$$

be the *first exit time* of $X$ from $D$. An exit from $D$ means that either $X$ enters $B_\varrho(x_+)$, that is, we have a transition to the other equilibrium, or $X$ leaves $B_R(0)$. But, as a consequence of our growth condition (1.16), the probability of the latter event does not contribute on the large deviations scale due to Proposition 1.4, as the following simple argument shows.

Let $s(\varepsilon, \mu) = sT^\varepsilon$ for some $s > 0$. Since $\tau_D = \tau_\varrho \wedge \sigma_R$, where $\sigma_R$ is the time of the diffusion's first exit from $B_R(0)$, Proposition 1.4 provides constants $R_1, \varepsilon_1 > 0$ s.t. for $R \geq R_1$, $\varepsilon \leq \varepsilon_1$,

$$\mathbb{P}_x(\tau_\varrho \geq s(\varepsilon, \mu)) \leq \mathbb{P}_x(\{\tau_\varrho \geq s(\varepsilon, \mu)\} \cap \{\sigma_R \geq s(\varepsilon, \mu)\}) + \mathbb{P}_x(\sigma_R < s(\varepsilon, \mu))$$

$$\leq \mathbb{P}_x(\tau_D \geq s(\varepsilon, \mu)) + C\eta^2 \frac{s(\varepsilon, \mu)}{\varepsilon} e^{-\eta R/\varepsilon} \qquad \text{for } \|x\| \leq \frac{R}{2}.$$

By the choice of $s(\varepsilon, \mu)$ and $T^\varepsilon = \exp(\frac{\mu}{\varepsilon})$, the right-hand side term in the last sum is of the order $\frac{1}{\varepsilon}\exp\frac{\mu - \eta R}{\varepsilon}$, that is, it can be assumed to be exponentially small of any exponential order required by choosing $R$ suitably large. Obviously, this holds uniformly with respect to $\mu$ on compact sets. This argument shows that the investigation of asymptotic properties of the laws of $\tau_\varrho$ may be replaced by a study of similar properties of $\tau_D$, with an error that may be chosen arbitrarily small by increasing $R$.

Similarly to the time homogeneous exit problem, we need a lemma which shows how to approximate the energy of a transition by the cost along particular trajectories which exit some neighborhood of $D$. This is of central importance to the estimation of the asymptotic law of $\tau_D$.

LEMMA 2.4. *Let $\vartheta > 0$ and $M$ a compact interval of $\mathbb{R}_+$. Then there exist $T_0 > 0$ and $\delta > 0$ with the following property:*

*For all $x \in D$ and $s \in M$, we can find a continuous path $\zeta^{x,s} \in C_{0T_0}$ starting in $\zeta_0^{x,s} = x$ and ending at some point of distance $d(\zeta_{T_0}^{x,s}, D) \geq \delta$ away from $D$ such that*

$$I_{0T_0}^s(\zeta^{x,s}) \leq e(s) + \vartheta \qquad \text{for all } s \in M.$$

PROOF. This proof extends arguments presented in Lemmas 5.7.18 and 5.7.19 in [4].

Fix $\vartheta > 0$, and let us decompose the domain $D$ into three different ones. Fixing $l > 0$, define a domain $\beta_l$ by $\beta_l = \{x \in D : \text{dist}(x, \chi) < l\}$. We recall that $\chi$ is the separation between $A_-$ and $A_+$. Then we define two closed sets $D_- = (D \setminus \beta_l) \cap A_-$ and $D_+ = (D \setminus \beta_l) \cap A_+$. We shall construct appropriate paths from points $y \in D$ to points a positive distance away from $D$ not exceeding the energy $e(s)$ by more than $\vartheta$ uniformly in $s \in M$ in four steps.

*Step* 1. Assume first that $y \in D_-$. For $l > 0$ small enough, we construct $\delta_1^l > 0$, $S_1^l > 0$ and a path $\psi_1^{s,y,l}$ defined on a time interval $[0, \tau_1^{s,y,l}]$ with $\tau_1^{s,y,l} \leq S_1^l$ for all $y \in D_-, s \in M$ and along which we exit a $\delta_1^l$-neighborhood of $D_-$ at cost at most $e(s) + \frac{2}{3}\vartheta$.

*Step* 1.1. In a first step we go from $y$ to a small neighborhood $B_l(x_-)$ of $x_-$, in time at most $T_1^l < \infty$, without cost.

We denote by $\varphi_1^{s,y,l}$ the trajectory starting at $\varphi_1^{s,y,l}(0) = y \in D_-$ of

$$\dot{\varphi}_1(t) = b(s, \varphi_1(t)),$$

and reaching $B_l(x_-)$ at time $\sigma_1^{y,s,l}$. Since $D_- \subset A_-$ and due to Assumption 2.1, $\sigma_1^{y,s,l}$ is finite. Moreover, since $b$ is locally Lipschitz, stability of solutions with respect to initial conditions and smooth changes of vector fields implies that there exist open neighborhoods $\mathcal{W}_y$ of $y$ and $\mathcal{W}_s$ of $s$ and $T_1^{s,y,l} > 0$ such that, for all $z \in \mathcal{W}_y$, $u \in \mathcal{W}_s$, $\sigma_1^{u,z,l} \leq T_1^{s,y,l}$. Recall that $D_-$ is compact. Therefore, we may find a finite cover of $D_- \times M$ by such sets, and, consequently, $T_1^l < \infty$ such that, for all $y \in D_-$ and $s \in M$, $\sigma_1^{s,y,l} \leq T_1^l$. Denote $z^{s,y,l} = \varphi_1^{s,y,l}(\sigma_1^{s,y,l})$.

*Step* 1.2. In a second step, we go from a small neighborhood $B_l(x_-)$ of $x_-$ to the equilibrium point $x_-$, in time at most 1, at cost at most $\frac{\vartheta}{3}$. In fact, by the continuity of the cost function, for $l$ small enough, $s \in M$, there exists a continuous path $\varphi_2^{s,y,l}$ of time length $\sigma_2^{s,y,l} \leq 1$ such that $\varphi_2^{s,y,l}(0) = z^{s,y,l}$, $\varphi_2^{s,y,l}(\sigma_2^{s,y,l}) = x_-$ and $I_{0\sigma_2^{s,y,l}}(\varphi_2^{s,y,l}) \leq \vartheta/3$.

*Step* 1.3. In a third step, we exit some $\delta$-neighborhood of $D_-$, starting from the equilibrium point $x_-$, in time at most $T_3 < \infty$, at cost at most $e(s) + \frac{\vartheta}{3}$ for $s \in M$.

By (2.4) and the continuity of the cost function for any $s \in M$, there exists $z_s \notin A_- \supset D_-$, $T_3^s < \infty$, some neighborhood $\mathcal{W}_s$ of $s$ and for $u \in \mathcal{W}_s$, we have $\varphi_3^u \in C_{0\sigma_3^u}$ such that $\varphi_3^u(0) = x_-$, $\varphi_3^u(\sigma_3^u) = z_s$, $\sigma_3^u \leq T_3^s$ and

$$\sup_{u \in \mathcal{W}_s} I_{0\sigma_3^u}^u(\varphi_3^u) \leq e(s) + \vartheta/3.$$

Use the compactness of $M$ to find a finite cover of $M$ by such neighborhoods, and, thus, some $T_3 < \infty$ such that all the statements hold with $\sigma_3^s \leq T_3$ for all $s \in M$. Finally, remark that the exit point is at least a distance $\tilde{\delta} = \inf_{i \in J} |z_i|$ away from the boundary of $D_-$, if $z_i, i \in J$, are the exit points corresponding to the finite cover.

In order to complete Step 1, we now define a path $\psi_1^{s,y,l} \in C_{0\tau_1^{s,y,l}}$ by concatenating $\varphi_1^{s,y,l}$, $\varphi_2^{s,y,l}$ and $\varphi_3^s$. This way, for small $l > 0$, we find $S_1^l > 0$ such that, for all $s \in M, y \in D_-$, we have $\tau_1^{s,y,l} \leq S_1^l$, $\psi_1^{s,y,l}(\tau_1^{s,y,l}) = y$, $\psi_1^{s,y,l}(\tau_1^{s,y,l}) \notin A_-$ and

$$I_{0\tau_1^{s,y,l}}^s(\psi_1^{s,y,l}) \leq e(s) + \tfrac{2}{3}\vartheta \qquad \text{for all } s \in M, y \in D_-.$$

At this point, we can encounter two cases. In the first case $\psi_1^{s,y,l}$ exits a $\delta_l$-neighborhood of $B_R(0)$. In this case we continue with Step 4. In the second case, $\psi_1^{s,y,l}$ exits $D_-$ into $\beta_l$, and we continue with Step 2.

*Step* 2. For $l$ small enough, we start in $y \in \beta_l$ to construct $S_2^l > 0$ and a path $\psi_2^{s,y,l}$ defined on a time interval $[0, \tau_2^{s,y,l}]$ with $\tau_2^{s,y,l} \leq S_2^l$ for all $y \in D_-, s \in M$ and along which we exit $\beta_l$ into the interior of $D_+$ at cost at most $\frac{\vartheta}{3}$.

In fact, due to the continuity of the cost function (see Lemma 1.6), there exists $l > 0$ small enough such that, for any $s \in M, y \in \overline{\beta_l}$, there exists $z_{s,y,l}$ in the interior of $D_+$, such that, $\psi_2^{s,y,l}(0) = y$ $\psi_2^{s,y,l}(\tau_2^{s,y,l}) = z_{s,y,l}$ and $I_{0\tau_2^{s,y,l}}^u(\psi_2^{s,y,l}) \leq \vartheta/3$. We may take $S_2^l = 1$.

*Step* 3. We start in $y \in D_+$ to construct $\delta_3^l > 0, S_3^l > 0$ and a path $\psi_3^{s,y,l}$ defined on a time interval $[0, \tau_3^{s,y,l}]$ with $\tau_3^{s,y,l} \leq S_3^l$ for all $y \in D_-, s \in M$ and along which we exit $D_+$ into $B_{\varrho-\delta_3^l}(x_+)$ at no cost.

Let $\delta_3^l = \varrho/2$. Since $D_+$ is compact and contained in the domain of attraction of $x_+$, stability of the solutions of the differential equation $\dot{\varphi}(t) = b(s, \varphi(t))$ with respect to the initial condition $y \in D_+$ and the parameter $s$ guarantees the existence of some time $S_3^l > 0$ such that the entrance time $\tau_3^{s,y,l}$ of $B_{\varrho/2}(x_+)$ by the solution starting in $y$ is bounded by $S_3^l$. Therefore, we may take $\psi_3^{s,y,l}$ to be defined by this solution restricted to the time interval before its entrance into $B_{\varrho/2}(x_+)$.

*Step* 4. For $l > 0$ small enough, we start in $x \in D_-$ and construct $T_0 > 0, \delta > 0$ and a path $\zeta^{s,x}$ defined on the time interval $[0, T_0]$, exiting a $\delta$-neighborhood of $D$ at cost at most $e(s) + \vartheta$ for all $s \in M$.

For $l$ small enough, take $T_0 = S_1^l + S_2^l + S_3^l$. We just have to concatenate paths constructed in the first three steps. Recall that $\psi_1^{s,x,l}$ passes through the equilibrium $x_-$ due to Step 1. In case $\psi_1^{s,x,l}$ exits a $\delta_1^l$-neighborhood of $B_R(0)$, just let the path spend enough time in $x_-$ without cost to obtain a path $\zeta^{s,x,l}$ defined on $[0, T_0]$, and take $\delta = \delta_1^l$. In the other case, we concatenate three paths constructed in Steps 1–3 to obtain a path defined on a subinterval of $[0, T_0]$ depending on $s, x, l$ and which exits a $\delta_3^l$-neighborhood of $D$. Recall from Step 1 that this path also passes through $x_-$. It remains to redefine the path by spending extra time at no cost in this equilibrium point to complete the proof. □

We now proceed to the estimation of uniform lower bounds for the asymptotic law of $\tau_D$. The uniformity has to be understood in the sense of Remark 1.9.

PROPOSITION 2.5. *Let $K$ be a compact subset of $D$:*

(a) *If $e(s) > \mu$, then*

$$\liminf_{\varepsilon \to 0} \varepsilon \log \inf_{x \in K} \mathbb{P}_x(\tau_D < sT^\varepsilon) \geq \mu - e(s),$$

*locally uniformly on* $\{(s, \mu) : \mu_* < \mu < \min(e(0), e(s)), 0 \leq s \leq 1\}$.

(b) *If $e(s) < \mu$, then*

$$\lim_{\varepsilon \to 0} \varepsilon \log \sup_{x \in K} \mathbb{P}_x(\tau_D \geq sT^\varepsilon) = -\infty,$$

*locally uniformly on* $\{(s, \mu) : e(s) < \mu < e(0), 0 \leq s \leq 1\}$.

PROOF. We choose a compact subset $L$ of $[0, 1]$ and a compact subset $M$ of $]\mu_*, e(0)[$, as well as some $\vartheta > 0$ such that

$$|e(s) - \mu| \geq \vartheta \qquad \forall (s, \mu) \in L \times M.$$

Later on we shall assume that $e(s) - \mu$ is uniformly positive respectively negative in order to prove (a) respectively (b).

In a first step, we apply Lemma 2.4 to approximate the energy function $e(s)$ by the rate function along a particular path, uniformly w.r.t. $s$. For the chosen $\vartheta$, it yields $T_0 > 0$ and $\delta > 0$, as well as continuous paths $\zeta^{x,s}$ indexed by $x \in D$ and $s \in [0, 1]$ ending a distance at least $\delta$ away from $D$ such that, for all $x \in D$ and $s \in [0, 1]$,

$$I_{0T_0}^s(\zeta^{x,s}) \leq e(s) + \frac{\vartheta}{4}.$$

In a second step, we use the Markov property to estimate the probability of exiting $D$ after time $sT^\varepsilon$ by a large product of exit probabilities after time intervals of length independent of $\varepsilon$ and $\mu$. Since for $\varepsilon > 0$, $\mu \in M$, the interval $[0, sT^\varepsilon]$ becomes arbitrarily large as $\varepsilon \to 0$, we introduce a splitting into intervals of length $\nu \geq T_0$ independent of $\varepsilon$ and $\mu$. For $k \in \mathbb{N}_0$, let $t_k = t_k(s, \varepsilon, \mu) := sT^\varepsilon - k\nu$. Then we have, for $k \in \mathbb{N}_0$ and $x \in D$,

$$\mathbb{P}_x(\tau_D \geq t_k) = \mathbb{E}_x(\mathbf{1}_{\{\tau_D \geq t_k\}} \mathbf{1}_{\tau_D \geq t_{k+1}})$$
$$= \mathbb{E}_x(\mathbf{1}_{\{\tau_D \geq t_{k+1}\}} \mathbb{E}[\mathbf{1}_{\{\tau_D \geq t_k\}} | \mathcal{F}_{t_{k+1}}])$$
$$\leq \mathbb{P}_x(\tau_D \geq t_{k+1}) \sup_{y \in D} \mathbb{P}_{y, t_{k+1}}(\tau_D \geq \nu).$$

Here $\mathbb{P}_{y,s}$ denotes the law of $X^{s,y}$, defined by the SDE

$$dX_t^{s,y} = b\left(\frac{s+t}{T^\varepsilon}, X_t^{s,y}\right) dt + \sqrt{\varepsilon}\, dW_t, \qquad t \geq 0,\ X_0^{s,y} = y \in \mathbb{R}^d.$$

On intervals $[0, \nu]$ it coincides with the law of the original process $X$ on $[s, s + \nu]$. Denoting $q_k(s, \varepsilon, \mu) := \sup_{y \in D} \mathbb{P}_{y, t_k}(\tau_D \geq \nu)$, an iteration of the latter argument yields

$$(2.7) \qquad \sup_{x \in K} \mathbb{P}_x(\tau_D \geq sT^\varepsilon) \leq \prod_{k=1}^{N(\varepsilon,\mu)} q_k(s, \varepsilon, \mu)$$

whenever $N(\varepsilon, \mu)\nu < sT^\varepsilon$. For the further estimation of the $q_k$, we apply some LDP to the product (2.7). This relies on the following idea. We choose $N(\varepsilon, \mu)$ of the order $\varepsilon T^\varepsilon$. Then the starting times $t_k$ appearing in the product belong to some neighborhood of $sT^\varepsilon$ that, compared to $T^\varepsilon$, shrinks to a point asymptotically. Consequently, the family of diffusions underlying the product is uniformly exponentially equivalent to the homogeneous diffusion governed by the drift $b(s, \cdot)$.

This will be done in the following third step. For $x \in D$, $s \in [0, 1]$, let

$$\Psi(x, s) := \left\{ \psi \in C_{0T_0} : \rho_{0T_0}(\psi, \zeta^{x,s}) < \frac{\delta}{2} \right\}$$

be the open $\delta/2$-neighborhood of the path chosen in the first step, and let

$$\Psi(x) := \bigcup_{s \in [0,1]} \Psi(x, s).$$

To apply our large deviations estimates in this situation, note first that conditions concerning $\tau_D$ translate into constraints for the trajectories of $X^\varepsilon$ as figuring in the preceding section: due to the definition of $\Psi(x,s)$, the choice $\nu \geq T_0$ and Lemma 2.4, we know that, for $y \in D, k \leq N(\varepsilon,\mu)$, if $X^{t_k,y}$ belongs to $\Psi(x)$, then for sure $X^{t_k,y}$ exits $D$ before time $\nu$. Keeping this in mind, we may apply Proposition 1.8 to the neighborhoods

$$S^{s,0}(\varepsilon,\mu) = [sT^\varepsilon - \nu N(\varepsilon,\mu), sT^\varepsilon + \nu]$$

of $sT^\varepsilon$. Each of the intervals $[t_k, t_k + \nu]$ is contained in $S^{s,0}(\varepsilon,\mu)$. As mentioned before, $N(\varepsilon,\mu)$ is chosen of the order $\varepsilon T^\varepsilon$, and this can be done uniformly w.r.t. $\mu \in M$. More precisely, we assume to have constants $0 < c_1 < c_2$ such that $c_1 \varepsilon T^\varepsilon \leq N(\varepsilon,\mu) \leq c_2 \varepsilon T^\varepsilon$. Then

$$\lim_{\varepsilon \to 0} \sup_{s \in [0,1], \mu \in M} \frac{\operatorname{diam} S^{s,0}(\varepsilon,\mu)}{T^\varepsilon} = 0,$$

and by the large deviations principle of Proposition 1.8, we obtain the lower bound

$$\liminf_{\varepsilon \to 0} \varepsilon \log \inf_{y \in K, \mu \in M, k \leq N(\varepsilon,\mu)} \mathbb{P}_{y,t_k}(\tau_D < \nu) \geq -\sup_{y \in K} \inf_{\psi \in \Psi(y)} I^s_{0T_0}(\psi)$$

$$\geq -\sup_{y \in K} I^s_{0T_0}(\zeta^{y,s}) \geq -e(s) - \frac{\vartheta}{4}.$$

We stress that this bound is uniform w.r.t. $s$ in the sense of Remark 1.9, so we can find $\varepsilon_0 > 0$ independent of $s$ such that, for $\varepsilon \leq \varepsilon_0$, $\mu \in M$ and $k \leq N(\varepsilon,\mu)$,

$$1 - q_k(s,\varepsilon,\mu) = \inf_{y \in D} \mathbb{P}_{y,t_k}(\tau_D < \nu)$$

$$\geq \inf_{y \in D, \mu \in M, j \leq N(\varepsilon,\mu)} \mathbb{P}_{y,t_j}(\tau_D < \nu) \geq \exp\left\{-\frac{1}{\varepsilon}\left(e(s) + \frac{\vartheta}{2}\right)\right\}.$$

From this, we obtain

$$\sup_{x \in K} \mathbb{P}_x(\tau_D \geq sT^\varepsilon) \leq \prod_{k=1}^{N(\varepsilon,\mu)} q_k(s,\varepsilon,\mu) \leq \left(1 - \exp\left\{-\frac{1}{\varepsilon}\left(e(s) + \frac{\vartheta}{2}\right)\right\}\right)^{N(\varepsilon,\mu)}$$

$$= \exp\left\{N(\varepsilon,\mu) \log\left(1 - \exp\left\{-\frac{1}{\varepsilon}\left(e(s) + \frac{\vartheta}{2}\right)\right\}\right)\right\} =: m(\varepsilon,\mu).$$

Since $\log(1-x) \leq -x$ for $0 \leq x < 1$, we have

$$m(\varepsilon,\mu) \leq \exp\left\{-c_1 \varepsilon \exp\left\{\frac{\mu}{\varepsilon} - \frac{1}{\varepsilon}\left(e(s) + \frac{\vartheta}{2}\right)\right\}\right\}.$$

In the fourth and last step, we exploit this bound of $m(\varepsilon,\mu)$ to obtain the claimed asymptotic bounds.

In order to prove (a), assume that $\mu < e(s)$ for $(s,\mu) \in L \times M$. Then the inner exponential approaches 0 on $L \times M$. Using the inequality $1 - e^{-x} \geq x \exp(-1)$ on $[0,1]$, we conclude that there exists $\varepsilon_1 \in (0, \varepsilon_0)$ such that, for all $\varepsilon \leq \varepsilon_1$ and $(s,\mu) \in L \times M$,

$$\varepsilon \log \inf_{x \in K} \mathbb{P}_x(\tau_D < sT^\varepsilon) \geq \varepsilon \log(1 - m(\varepsilon, \mu))$$

$$\geq \varepsilon \log\left(\varepsilon c_1 \exp(-1) \exp\left\{\frac{1}{\varepsilon}\left(\mu - e(s) - \frac{\vartheta}{2}\right)\right\}\right)$$

$$= -\varepsilon + \varepsilon \log c_1 + \varepsilon \log \varepsilon + \mu - e(s) - \frac{\vartheta}{2}$$

$$\geq \mu - e(s) - \vartheta.$$

For (b), assume $\mu > e(s)$ on $L \times M$. Then

$$\varepsilon \log \sup_{x \in K} \mathbb{P}_x(\tau_D \geq sT^\varepsilon) \leq \varepsilon \log m(\varepsilon, \mu)$$

$$\leq -c_1 \varepsilon \exp\left\{-\frac{1}{\varepsilon}\left(\mu - e(s) - \frac{\vartheta}{2}\right)\right\} \xrightarrow[\varepsilon \to 0]{} -\infty. \quad \square$$

As a consequence of these large deviations type results on the asymptotic distribution of $\tau_D$ and the remarks preceding the statement of Lemma 2.4 and Proposition 2.5, we get the following asymptotics for the transition time of the diffusion.

PROPOSITION 2.6. *Let $x \in A_-$. There exists $h_0 > 0$ such that*

(2.8) $$\liminf_{\varepsilon \to 0} \varepsilon \log \mathbb{P}_x(\tau_\varrho \leq (a_\mu - h)T^\varepsilon) \geq \mu - e(a_\mu - h),$$

(2.9) $$\lim_{\varepsilon \to 0} \varepsilon \log \mathbb{P}_x(\tau_\varrho \geq (\alpha_\mu + h)T^\varepsilon) = -\infty,$$

*for $h \leq h_0$. Moreover, these convergence statements hold uniformly w.r.t. $x$ on compact subsets of $D$ and w.r.t. $\mu$ on compact subsets of $]\mu_*, e(0)[$.*

PROOF. As the estimation based on Proposition 1.4 at the beginning of the section shows, we may derive the required estimates for $\tau_D$ instead of $\tau_\varrho$, if $R$ is chosen large enough.

Let $M$ be a compact subset of $]\mu_*, e(0)[$. Then $0 < a_\mu < 1$ for $\mu \in M$ which yields the existence of $h_0 > 0$ such that the compact set $L_h := \{a_\mu - h : \mu \in M\}$ is contained in $]0,1[$ for $h \leq h_0$. Moreover, we have $e(s) > \mu$ for $0 < s < a_\mu$ due to the assumptions on $e$, uniformly w.r.t. $(s, \mu) \in L_h \times M$ by the continuity of $e$. Hence, by Proposition 2.5(a),

$$\liminf_{\varepsilon \to 0} \varepsilon \log \inf_{x \in K} \mathbb{P}_x(\tau_D \leq sT^\varepsilon) \geq \mu - e(s),$$

uniformly on $L_h \times M$ for all $h \leq h_0$. By setting $s = a_\mu - h$, we obtain the first asymptotic inequality. The second one follows in a completely analogous way from Proposition 2.5(b) since $\alpha_\mu = a_\mu$ and $e(a_\mu + h) < \mu$ for small enough $h$. □

2.2. *Upper bound for the exponential exit rate.* Let us next derive upper bounds for the exponential exit rate which resemble the lower bounds just obtained. We need an extension of a result obtained by Freidlin and Wentzell (Lemma 5.4 in [18]).

LEMMA 2.7. *Let $K$ be a compact subset of $A_- \setminus \{x_-\}$. There exist $T_0 > 0$ and $c > 0$ such that, for all $T \geq T_0$, $s \in [0,1]$ and for each $\varphi \in C_{0T}$ taking its values in $K$, we have*

$$I^s_{0T}(\varphi) \geq c(T - T_0).$$

PROOF. Let $\phi_{s,x}$ be the solution of the differential equation

$$\dot\phi_{s,x}(t) = b(s, \phi_{s,x}(t)), \qquad \phi_{s,x}(0) = x \in K.$$

Let $\tau(s,x)$ be the first exit time of the path $\phi_{s,x}$ from the domain $K$. Since $A_-$ is the domain of attraction of $x_-$ and since $K$ is a compact subset of $A_- \setminus \{x_-\}$, we obtain $\tau(s,x) < \infty$ for all $x \in K$.

The function $\tau(s,x)$ is upper semi-continuous with respect to the variables $s$ and $x$ (due to the continuous dependence of $\phi_{s,x}$ on $s$ and $x$). Hence, the maximal value $T_1 := \sup_{s \in [0,1], \, x \in K} \tau(a, x)$ is attained.

Let $T_0 = T_1 + 1$, and consider all functions $\varphi \in C_{0T_0}$ with values in $K$. This set of functions is closed with respect to the maximum norm. Since there is no solution of the ordinary differential equation in this set of functions, the functional $I^s_{0T_0}$ reaches a strictly positive minimum on this set which is uniform in $s$. Let us denote it by $m$. By the additivity of the functional $I^s_{0T}$, we obtain, for $T \geq T_0$ and $\varphi \in C_{0T}$ with values in $K$,

$$I^s_{0T}(\varphi) \geq m \left\lfloor \frac{T}{T_0} \right\rfloor \geq m \left( \frac{T}{T_0} - 1 \right) = c(T - T_0),$$

with $c = \frac{m}{T_0}$. □

Let us recall the subject of interest of this subsection:

$$\tau_\varrho = \inf\{t \geq 0 : X^\varepsilon_t \in B_\varrho(x_+)\},$$

the hitting time of a small neighborhood of the equilibrium point $x_+$. First we shall consider upper bounds for the law of this time in some window of length $\beta T^\varepsilon$ where $\beta$ is sufficiently small. The important feature of the following statement is that $\beta$ is independent of $s$, while the uniformity of the bound again has to be understood in the sense of Remark 1.9.

PROPOSITION 2.8. *For all $\vartheta > 0$, there exist $\beta > 0$, $\eta > 0$ such that, for all $s \in [0,1]$,*

$$\limsup_{\varepsilon \to 0} \varepsilon \log \sup_{x \in B_\eta(x_-)} \mathbb{P}_x(sT^\varepsilon \leq \tau_\varrho \leq (s+\beta)T^\varepsilon) \leq \mu - e(s) + \vartheta.$$

*This bound holds locally uniformly w.r.t. $\mu \in ]\mu_*, e(0)[$ and uniformly w.r.t. $s \in [0,1]$.*

PROOF. Let $M$ be a compact subset of $]\mu_*, e(0)[$, and fix $\vartheta > 0$. We first introduce some parameter dependent domains, the exit times of which will prove to be suitable for estimating the probability that $\tau_\varrho$ is in a certain time window.

For this purpose, we define for $\delta > 0$ and $s \in [0,1]$ an open domain

$$D(\delta, s) := \left\{ y \in \mathbb{R}^d : V^s(x_-, y) < \mu^* + \frac{1}{1+\delta}, \operatorname{dist}(y, A_+) > \delta \right\},$$

and we let $D = D(\delta) = \bigcup_{s \in [0,1]} D(\delta, s)$. Then $D$ is relatively compact in $A_-$, $\operatorname{dist}(y, A_+) > \delta$ for all $y \in D(\delta)$, and a transition to a $\varrho$-neighborhood of $x_+$ certainly requires an exit from $D(\delta)$. The boundary of $D(\delta)$ consists of two hyper surfaces, one of which carries an energy strictly greater than $\mu^*$ and thus greater than $e(s)$ for all $s \in [0,1]$. The minimal energy is therefore attained on the other component of distance $\delta$ from $A_+$ which approaches $\chi = \partial A_-$ as $\delta \to 0$. Thus, by the joint continuity of the quasi-potential, we can choose $\delta_0 > 0$ and $\eta > 0$ such that, for $\delta \leq \delta_0$ and $s \in [0,1]$,

(2.10)
$$e(s) = \inf_{z \in \chi} V^s(x_-, z) \geq \inf_{z \in \partial D(\delta)} V^s(x_-, z)$$

$$\geq \inf_{y \in B_\eta(x_-)} \inf_{z \in \partial D(\delta)} V^s(y, z) \geq e(s) - \frac{\vartheta}{4}.$$

Let $\tau_D$ be the first exit time of $X^\varepsilon$ from $D$. For $s \in [0,1]$ and $\beta > 0$, we introduce a covering of the interval of interest $[sT^\varepsilon, (s+\beta)T^\varepsilon]$ into $N = N(\beta, \varepsilon, \mu)$ intervals of fixed length $\nu$, that is, $\nu$ is independent of $\varepsilon$, $\mu$, $s$ and $\beta$. We will have to assume that $\nu$ is sufficiently large, which will be made precise later on. Thus, we have $N\nu \geq \beta T^\varepsilon$, and we can and do assume that $N \leq \beta T^\varepsilon$. For $k \in \mathbb{Z}$, $k \geq -1$, let

$$t_k = t_k(s, \varepsilon, \mu) := sT^\varepsilon + k\nu$$

be the starting points of these intervals. We consider $t_{-1}$ since we need some information about the past in order to ensure the diffusion to start in a neighborhood of the equilibrium $x_-$. Then for $x \in B_\eta(x_-)$, we get the desired estimation of probabilities of exit windows for $\tau_\varrho$ by those with respect to $\tau_D$:

$$\mathbb{P}_x(sT^\varepsilon \leq \tau_\varrho \leq (s+\beta)T^\varepsilon) \leq \sum_{k=0}^{N} \mathbb{P}_x(t_k \leq \tau_D \leq t_{k+1}).$$

In a second step we will fix $k \geq 0$ and estimate the probability of a first exit from $D$ during each of the intervals $[t_k, t_{k+1}]$ separately. Here the difficulty is that we do not have any information on the location at time $t_k$. We therefore condition on whether or not $X^\varepsilon$ has entered the neighborhood $B_\eta(x_-)$ in the previous time interval. For that purpose, let

$$\sigma_k := \inf\{t \geq t_k \vee 0 : X_t^\varepsilon \in B_\eta(x_-)\}, \qquad k \geq -1.$$

Then for $k \geq 0$,

$$\begin{aligned}(2.11) \quad \mathbb{P}_x(t_k \leq \tau_D \leq t_{k+1}) \leq\ & \mathbb{P}_x(t_k \leq \tau_D \leq t_{k+1} | \sigma_{k-1} \leq t_k) \\ & + \mathbb{P}_x(\tau_D \wedge \sigma_{k-1} \geq t_k).\end{aligned}$$

In the next step we shall estimate the second term on the right-hand side of (2.11). Let $K = K(\delta, \eta) = \overline{D(\delta)} \setminus B_\eta(x_-)$. Then $K$ is compact, and by the Markov property, we have

$$\mathbb{P}_x(\tau_D \wedge \sigma_{k-1} \geq t_k) \leq \sup_{y \in K} \mathbb{P}_{y, t_{k-1}}(\tau_D \wedge \sigma_1 \geq \nu),$$

where $\mathbb{P}_{y,t}$ is as defined in the previous section. Now we wish to further estimate this exit probability using large deviations methods. The neighborhoods

$$S^{s,\beta}(\varepsilon, \mu) = [sT^\varepsilon - \nu, sT^\varepsilon + \nu N(\beta, \varepsilon, \mu)]$$

of $sT^\varepsilon$ contain each interval $[t_k, t_{k+1}]$, $-1 \leq k \leq N(\beta, \varepsilon, \mu)$, and they satisfy

$$\limsup_{\varepsilon \to 0} \sup_{\mu \in M, s \in [0,1]} \frac{\mathrm{diam}(S^{s,\beta}(\varepsilon, \mu))}{T^\varepsilon} \leq \beta.$$

Hence, by the uniform LDP of Proposition 1.8, applied to the closed set

$$\Phi_K(\delta, \eta) = \{\varphi \in C_{0,\nu} : \varphi_t \in K(\delta, \eta) \text{ for all } t \in [0, \nu]\},$$

we obtain the upper bound

$$\begin{aligned}(2.12) \quad & \limsup_{\varepsilon \to 0} \varepsilon \log \sup_{y \in K, \mu \in M, k \leq N} \mathbb{P}_{y, t_{k-1}}(\tau_D \wedge \sigma_1 \geq \nu) \\ & \leq \limsup_{\varepsilon \to 0} \varepsilon \log \sup_{y \in K, \mu \in M, t \in S^{s,\beta}(\varepsilon, \mu)} \mathbb{P}_{y,t}(X^\varepsilon \in \Phi_K(\delta, \eta)) \\ & \leq - \inf_{y \in K} \inf_{\varphi \in \Phi_K(\delta, \eta)^{\gamma_0(\beta)}} I_{0,\nu}^s(\varphi),\end{aligned}$$

where $\gamma_0(\beta) = \beta \delta(\Phi_K(\delta, \eta))$ is the "blowup-factor" induced by the diameter $\beta$. Since $\gamma_0(\beta) \to 0$ as $\beta \to 0$, we can find $\beta_0 > 0$ such that, for $\beta \leq \beta_0$,

$$\Phi_K(\delta, \eta)^{\gamma_0(\beta)} \subset \Phi_K\left(\frac{\delta}{2}, \frac{\eta}{2}\right),$$

which amounts to saying that, instead of blowing up the set of paths, we consider the slightly enlarged domain $K(\frac{\delta}{2}, \frac{\eta}{2})$. Thus,

$$-\inf_{y \in K} \inf_{\varphi \in \Phi_K(\delta,\eta)^{\gamma_0(\beta)}} I^s_{0,\nu}(\varphi) \leq -\inf_{y \in K} \inf_{\varphi \in \Phi_K(\frac{\delta}{2},\frac{\eta}{2})} I^s_{0,\nu}(\varphi).$$

By Lemma 2.7, the latter expression, and therefore the r.h.s. of (2.12), approaches $-\infty$ as $\nu \to \infty$, uniformly w.r.t. $s \in [0,1]$. So the second term in the decomposition of $\mathbb{P}_x(t_k \leq \tau_D \leq t_{k+1})$ can be neglected since it becomes exponentially small of any desired order by choosing $\nu$ suitably large.

In the next and most difficult step, we treat the first term on the right-hand side of (2.11). It is given by the probability that, while $X^\varepsilon$ is in $B_\eta(x_-)$ at time $\sigma_{k-1}$, it exits within a time interval of length $t_{k+1} - \sigma_{k-1} \leq 2\nu$. Hence, by the strong Markov property,

$$\mathbb{P}_x(t_k \leq \tau_D \leq t_{k+1}|\sigma_{k-1} \leq t_k) \leq \sup_{t_{k-1} \leq t \leq t_k, y \in B_\eta(x_-)} \mathbb{P}_{y,t}(\tau_D \leq 2\nu).$$

Applying the uniform LDP to the closed set

$$F_D(\delta) := \{\varphi \in C_{0,2\nu} : \varphi_0 \in D(\delta), \varphi_{t_0} \notin D(\delta) \text{ for some } t_0 \leq 2\nu\}$$

yields the upper bound

(2.13)
$$\limsup_{\varepsilon \to 0} \varepsilon \log \sup_{y \in B_\eta(x_-), \mu \in M, t \in S^{s,\beta}(\varepsilon,\mu)} \mathbb{P}_{y,t}(\tau_D \leq 2\nu)$$
$$\leq -\inf_{y \in B_\eta(x_-)} \inf_{\varphi \in F_D(\delta)^{\gamma_0(\beta)}} I^s_{0,2\nu}(\varphi),$$

where $\gamma_0(\beta) = 2\beta\delta(F_D(\delta))$. By the same reasoning as before, we can replace the blow-up of the path sets $F_D(\delta)$ by an enlargement of the domain $D(\delta)$. We find $\beta_1 > 0$ such that, for $\beta \leq \beta_1$

$$-\inf_{y \in B_\eta(x_-)} \inf_{\varphi \in F_D(\delta)^{\gamma_0(\beta)}} I^s_{0,2\nu}(\varphi) \leq -\inf_{y \in B_\eta(x_-)} \inf_{\varphi \in F_D(\delta/2)} I^s_{0,2\nu}(\varphi)$$
$$\leq -\inf_{y \in B_\eta(x_-)} \inf_{z \in \partial D(\delta/2)} V^s(y,z).$$

Now we apply (2.10) and recall the uniformity of the LDP w.r.t. $s$. We find $\varepsilon_0 > 0$ such that we have, for $\varepsilon \leq \varepsilon_0$, $s \in [0,1]$, $\mu \in M$ and $\beta \leq \beta_1$,

(2.14)
$$\varepsilon \log \sup_{y \in B_\eta(x_-), t \in S^{s,\beta}(\varepsilon,\mu)} \mathbb{P}_{y,t}(\tau_D \leq 2\nu)$$
$$\leq -\inf_{y \in B_\eta(x_-)} \inf_{z \in \partial D(\delta/2)} V^s(y,z) + \frac{\vartheta}{4}$$
$$\leq -e(s) + \frac{\vartheta}{2}.$$

We finally summarize our findings. We conclude that there exists $\varepsilon_1 > 0$ such that, for $\varepsilon \leq \varepsilon_1$, $\mu \in M$ and $s \in [0,1]$, we have

$$\varepsilon \log \sup_{x \in B_\eta(x_-)} \mathbb{P}_x(sT^\varepsilon \leq \tau_\varrho \leq (s+\beta)T^\varepsilon)$$

$$\leq \varepsilon \log \left\{ \sum_{k=0}^{N(\beta,\varepsilon,\mu)} \sup_{x \in B_\eta(x_-)} \mathbb{P}_x(t_k \leq \tau_D \leq t_{k+1} | \sigma_{k-1} \leq t_k) \right\} + \frac{\vartheta}{4}$$

$$\leq \varepsilon \log \left\{ \beta T^\varepsilon \exp\left(-\frac{1}{\varepsilon}\left[e(s) - \frac{\vartheta}{2}\right]\right)\right\} + \frac{\vartheta}{4} = \varepsilon \log \beta + \mu - e(s) + \frac{3}{4}\vartheta$$

$$\leq \mu - e(s) + \vartheta.$$

This completes the proof. □

REMARK 2.9. If we stay away from $s = 0$, in the statement of Proposition 2.8 the radius of the starting domain $B_\eta(x_-)$ can be chosen independently of the parameter $\vartheta$. It may then be brought into the following somewhat different form.

PROPOSITION 2.10. *Let $L$ and $M$ be compact subsets of $]0,1]$ respectively $]\mu_*, e(0)[$. Let $\eta > 0$ be small enough such that $B_\eta(x_-)$ belongs to the domain*

$$\{y \in \mathbb{R}^d : V^s(x_-, y) < \mu^* \text{ for all } s \in L\}.$$

*Then, for all $\vartheta > 0$, there exists some $\beta > 0$ such that we have*

$$\limsup_{\varepsilon \to 0} \varepsilon \log \sup_{x \in B_\eta(x_-)} \mathbb{P}(sT^\varepsilon \leq \tau_\varrho \leq (s+\beta)T^\varepsilon) \leq \mu - e(s) + \vartheta,$$

*uniformly w.r.t. $s \in L$ and $\mu \in M$.*

PROOF. To prove Proposition 2.10, one has to modify slightly the preceding proof. Instead of just $\eta$, one has to choose two different parameters: $\eta_0$ for the definition of the starting domain $D$ and some $\eta_1$ for the description of the location of the diffusion at time $t_k$, that is, for the definition of the stopping times $\sigma_k$. □

In the following proposition, we derive the upper bound for the asymptotic law of transition times, corresponding to the lower bound obtained in Proposition 2.6.

PROPOSITION 2.11. *Let $\mu < e(0)$, and recall from (2.5) the definition $a_\mu = \inf\{t \geq 0 : e(t) \leq \mu\}$. There exist $\gamma > 0$ and $h_0 > 0$ such that, for all $h \leq h_0$,*

$$(2.15) \quad \limsup_{\varepsilon \to 0} \varepsilon \log \sup_{x \in B_\gamma(x_-)} \mathbb{P}_x(\tau_\varrho \leq (a_\mu - h)T^\varepsilon) \leq \mu - e(a_\mu - h).$$

*This bound is uniform w.r.t. $\mu$ on compact subsets of $]\mu_*, e(0)[$.*

PROOF. Let $M$ be a compact subset of $]\mu_*, e(0)[$. To choose $h_0$, we use our assumptions on the geometry of the energy function $e$. Recall Assumption 2.2 according to which $e$ is strictly monotonous in the open intervals between the extrema $]\mu_*, \mu^*[$. It implies that $e$ is monotonically decreasing on the interval $[a_{e(0)}, a_\mu]$ for any $\mu \in M$. By choice of $M$, we further have $a_{e(0)} < \inf_{\mu \in M} a_\mu$. Now choose $h_0$ such that $\inf_{\mu \in M} a_\mu - h_0 > a_{e(0)}$. Then we have, for $h \leq h_0$,

(2.16) $$\inf_{\mu \in M} a_\mu - h > 0,$$

(2.17) $$e(0) > \sup_{\mu \in M, h \leq h_0} e(a_\mu - h),$$

(2.18) $$e(s) \geq e(a_\mu - h) \quad \text{for all } s \leq a_\mu - h.$$

To see (2.18), note that, for $0 \leq s \leq a_{e(0)}$, by definition of $a_{e(0)}$, the inequality $e(s) \geq e(0) > e(a_\mu - h)$ holds, while, for $a_{e(0)} \leq s \leq a_\mu - h$ by monotonicity, $e(s) \geq e(a_\mu - h)$.

Next fix $h \leq h_0$. For $\mu \in M$, let $\Lambda_0 = \Lambda_0(\mu) = 0$, and $\Lambda_1(\mu) \leq \inf_{\mu \in M}(a_\mu - h)T^\varepsilon$. For $N \in \mathbb{N}^*$, we set $\Lambda_i(\mu) = \Lambda_1 + \frac{i-1}{N-1}\left((a_\mu - h)T^\varepsilon - \Lambda_1\right)$, $2 \leq i \leq N$, thus splitting the time interval $[0, (a_\mu - h)T^\varepsilon]$ into the $N$ intervals $[\Lambda_i(\mu), \Lambda_{i+1}(\mu)]$, $0 \leq i \leq N-1$. Then for $\gamma > 0$, $x \in B_\gamma(x_-)$,

$$\mathbb{P}_x(\tau_\varrho \leq (a_\mu - h)T^\varepsilon) \leq \sum_{i=0}^{N-1} \mathbb{P}_x(\tau_\varrho \in [\Lambda_i(\mu), \Lambda_{i+1}(\mu)]),$$

which implies

$$\limsup_{\varepsilon \to 0} \varepsilon \log \sup_{x \in B_\gamma(x_-)} \mathbb{P}_x(\tau_\varrho \leq (a_\mu - h)T^\varepsilon)$$
$$\leq \max_{0 \leq i \leq N-1} \limsup_{\varepsilon \to 0} \varepsilon \log \sup_{x \in B_\gamma(x_-)} \mathbb{P}_x(\tau_\varrho \in [\Lambda_i(\mu), \Lambda_{i+1}(\mu)]).$$

Fix $\vartheta > 0$ such that, for $h \leq h_0, \mu \in M$, we have $e(0) \geq e(a_\mu - h) + \vartheta$. This is guaranteed by (2.17). We shall show that

$$\limsup_{\varepsilon \to 0} \varepsilon \log \sup_{x \in B_\gamma(x_-)} \mathbb{P}_x(\tau_\varrho \in [\Lambda_i(\mu), \Lambda_{i+1}(\mu)]) \leq \mu - e(a_\mu - h) + \vartheta,$$

uniformly in $0 \leq i \leq N-1$ and $\mu \in M$.

Let us treat the estimation of the first term separately from the others. In fact, by Proposition 2.8, setting $s = 0$, $\beta = \Lambda_1/T^\varepsilon$, we may choose $\Lambda, \varepsilon_0 > 0$ and $\gamma_0 > 0$ such that, for $\Lambda_1 \leq \Lambda T^\varepsilon, \varepsilon \leq \varepsilon_0, \gamma \leq \gamma_0, \mu \in M$, the inequality

$$\varepsilon \log \sup_{x \in B_\gamma(x_-)} \mathbb{P}_x(\tau_\varrho \in [\Lambda_0(\mu), \Lambda_1(\mu)]) \leq \mu - e(0) + \vartheta$$

holds. Now we use the inequality $e(0) \geq e(a_\mu - h) + \vartheta$, valid for all $\mu \in M$. Hence, there exists $\Lambda > 0$, $\varepsilon_0 > 0$ and $\gamma_0 > 0$ such that, for $\Lambda_1 \leq \Lambda T^\varepsilon$, $\varepsilon \leq \varepsilon_0$, $\gamma \leq \gamma_0, \mu \in M$,

$$\varepsilon \log \sup_{x \in B_\gamma(x_-)} \mathbb{P}_x(\tau_\varrho \in [\Lambda_0(\mu), \Lambda_1(\mu)]) \leq \mu - e(a_\mu - h).$$

Let us next estimate the contributions for the intervals $[\Lambda_i(\mu), \Lambda_{i+1}(\mu)]$ with $i \geq 1$. We use Proposition 2.8, this time with $s = \Lambda_i(\mu)/T^\varepsilon$, $\beta = \frac{1}{N-1} \sup_{\mu \in M} a_\mu$. By the definition of $a_\mu$, we get $e(s) > e(a_\mu)$ for all $s < a_\mu$. By (2.18), we have $e(s) = e(\Lambda_i(\mu)/T^\varepsilon) \geq e(a_\mu - h)$. By Remark 2.9,

$$\limsup_{\varepsilon \to 0} \varepsilon \log \sup_{x \in B_\gamma(x_-)} \mathbb{P}_x(\tau_\varrho \in [\Lambda_i(\mu), \Lambda_{i+1}(\mu)]) \leq \mu - e(a_\mu - h) + \vartheta,$$

uniformly w.r.t. $1 \leq i \leq N$ and $\mu \in M$. Letting $\vartheta$ tend to 0, which implies that $N$ tends to infinity and $\Lambda_1$ tends to zero, we obtain the desired upper bound for the exponential exit rate. □

**3. Stochastic resonance.** Given the results of the previous section on the asymptotics of exit times which are uniform in the scale parameter $\mu$, we are now in a position to reconsider the problem of finding a satisfactory probabilistic notion of stochastic resonance that does not suffer from the lack of robustness defect of physical notions, such as spectral power amplification. We continue to study the SDE

$$dX_t^\varepsilon = b\left(\frac{t}{T^\varepsilon}, X_t^\varepsilon\right) dt + \sqrt{\varepsilon}\, dW_t, \qquad t \geq 0, X_0^\varepsilon = x_0 \in \mathbb{R}^d,$$

introduced before, thereby recalling that the drift term $b$ satisfies the local Lipschitz conditions (1.15) and (1.14) in space and time, as well as the growth condition (1.16). Moreover, $b$ is assumed to be one-periodic in time such that $T^\varepsilon$ is the period of the deterministic input of the randomly perturbed dynamical system described by $X^\varepsilon$.

In typical applications, $b = -\nabla_x U$ is given by the (spatial) gradient of some time periodic double-well potential $U$ (see [13]). The potential possesses at all times two local minima well separated by a barrier. The depth of the wells and the roles of being the deep and shallow one change periodically. The diffusion $X^\varepsilon$ then roughly describes the motion of a Brownian particle of intensity $\varepsilon$ in a double-well landscape. Its attempts to get close to the energetically most favorable deep position in the landscape makes it move along random trajectories which exhibit randomly periodic hopping between the wells. The average time the trajectories need to leave a potential well of depth $\frac{v}{2}$ being given by the Kramers–Eyring law $T^\varepsilon = \exp(\frac{v}{\varepsilon})$ motivates our choice of time scales $T^\varepsilon = \exp(\frac{\mu}{\varepsilon})$ and also our convention to measure time scales in energy units $\mu$.

The problem of stochastic resonance consists of characterizing the *optimal tuning* of the noise, that is, the best relation between the noise amplitude $\varepsilon$ and the input period $T^\varepsilon$—or, in our units the energy parameter $\mu$—of the deterministic system which makes the diffusion trajectories look as periodic as possible. Of course, the optimality criterion must be based upon a *quality measure* for periodicity in random trajectories.

In this section we shall develop a measure of quality based on the transition probabilities investigated in Section 2 and with respect to this measure for fixed small $\varepsilon$ (in the small noise limit $\varepsilon \to 0$), exhibit a resonance energy $\mu_0(\varepsilon)$ for which the diffusion trajectories follow the periodic forcing of the system at intensity $\varepsilon$ in an optimal way. We shall in fact study the problem in a more general situation which includes the double-well potential gradient case as an important example, and draws its intuition from it. The deterministic system $\dot\varphi_s(t) = b(s, \varphi_s(t)), t \geq 0$, has to satisfy Assumption 2.1, that is, it possesses two well separated domains of attraction, the common boundary of which is time invariant. In the first subsection we shall describe the *resonance interval*, that is, the set of all parameter values $\mu$ for which in the small noise limit trivial behavior, that is, either constant or continuously jumping trajectories, are excluded. The second subsection shows that a quality measure of periodic tuning is given by the exponential rate at which the first transition from one domain of attraction to the other one happens within a fixed time window around $a_\mu T^\varepsilon$. This quality measure is robust, as demonstrated in the last subsection: in the small noise limit the diffusion and its reduced model, a Markov chain jumping between the domains of attraction reduced to the equilibrium points, have the same resonance pattern.

3.1. *Resonance interval.* According to Freidlin [5], quasi-periodic hopping behavior of the trajectories of our diffusion in the small noise limit of course requires that the energies required to leave the domains of attraction of the two equilibria switch their order periodically: if $e_\pm$ denotes the energy needed to leave $A_\pm$, then $e_+$ needs to be bigger than $e_-$ during part of one period, and vice versa for the rest. We assume that $e_\pm$ both satisfy Assumption 2.2 and associate to each of these functions the transition time

$$a_\mu^\pm(s) = \inf\{t \geq s : e_\pm(t) \leq \mu\}.$$

The time scales $\mu$ for which relevant behavior of the system is expected clearly belong to the intervals

$$I_i = \left]\inf_{t \geq 0} e_i(t), \sup_{t \geq 0} e_i(t)\right[, \qquad i \in \{-, +\}.$$

Our aim being the observation of periodicity, we have to make sure that the process can travel back and forth between the domains of attraction on the time scales considered, but not instantaneously. So, on the one hand, in these scales it should not get stuck in one of the domains. On the other hand, they should not allow for chaotic behavior, that is, immediate re-bouncing after leaving a domain has to be avoided.

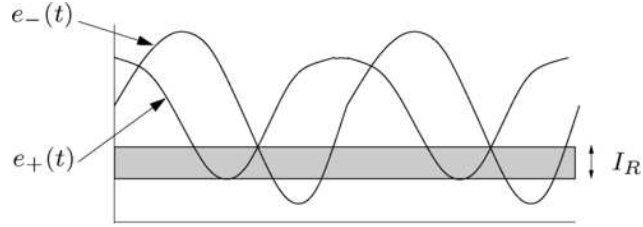

Fig. 3. *Resonance interval.*

To make these conditions mathematically precise, recall that transitions from $A_i$ to $A_{-i}$ become possible as soon as the energy $e_i$ needed to exit from domain $i$ falls below $\mu$ which represents the available energy. Not to get stuck in one of $A_\pm$, we therefore have to guarantee

$$\mu > \max_{i=\pm} \inf_{t\geq 0} e_i(t).$$

To avoid immediate re-bouncing, we have to assure that the diffusion cannot leave $A_{-i}$ at the moment it reaches it, coming from $A_i$. Suppose we consider the dynamics after time $s \geq 0$, and the diffusion is near $i$ at that time. Its first transition to $A_{-i}$ occurs at time $a^i_\mu(s)T^\varepsilon$, where $a^i_\mu(s)$ is the first time in the original scale at which $e_i$ falls below $\mu$ after $s$. Provided $e_{-i}(a^i_\mu(s))$ is bigger than $\mu$, it stays there for at least a little while. This is equivalent to saying that, for all $s \geq 0$, there exists $\delta > 0$ such that on $[a^i_\mu(s), a^i_\mu(s)+\delta]$ we have $\mu < e_{-i}$. Since by definition for $t$ shortly after $a^i_\mu(s)$, we always have $e_i(t) \leq \mu$, our condition may be paraphrased by the following: for all $s \geq 0$, there exists $\delta > 0$ such that on $[a^i_\mu(s), a^i_\mu(s)+\delta]$ we have $\mu < \max_{i=\pm} e_i$. This, in turn, is more elegantly expressed by

$$\mu < \inf_{t\geq 0} \max_{i=\pm} e_i(t).$$

Our search for a set of scales $\mu$ for which the diffusion exhibits nontrivial transition behavior may be summarized in the following definition. The interval

$$I_R = \left]\max_{i=\pm}\inf_{t\geq 0} e_i(t), \inf_{t\geq 0}\max_{i=\pm} e_i(t)\right[$$

is called *resonance interval* (see Figure 3).

In this interval, for small $\varepsilon$, we have to look for an optimal energy scale $\mu(\varepsilon)$ in the following subsection. See [9] and [8] for the definition of the corresponding interval in the one-dimensional case and in the case of two state Markov chains. In Freidlin's [5] terms, stochastic resonance in the sense of quasi-deterministic periodic motion is given if the parameter $\mu$ exceeds the lower boundary of our resonance interval.

Let us consider the potential gradient case. Assume that $b(t,x) = -\nabla_x U(t,x)$, $t \geq 0$, $x \in \mathbb{R}^d$, where $U$ is some time periodic double-well potential with time invariant local minima $x_\pm$ and separatrix. Then $A_-$ and $A_+$ represent the two wells of the potential, $\chi$ the separatrix. The energy $e_\pm$ is, in fact, the energy some Brownian particle needs to cross $\chi$. Freidlin and Wentzell [6] give the link between this energy and the depth of the well.

LEMMA 3.1. *If $D_\pm(t) = \inf_{y \in \chi} U(t,y) - U(t,x_\pm)$ denote the depths of the wells, then $e_\pm(t) = 2D_\pm(t)$ for all $t \geq 0$.*

This link is the origin of the name "quasipotential." The minimal energy $e$ is reached by some path which intersects the level sets of the potential with orthogonal tangents. This path satisfies an equation of the form

$$\dot\varphi_s = \nabla_x U(t, \varphi_s), \qquad s \in (-\infty, T),\ \varphi_T \in \chi.$$

The resonance interval is given by $I_R = \,]\max_{i=\pm} \inf_{t \geq 0} 2D_i(t), \inf_{t \geq 0} \max_{i=\pm 1} 2D_i(t)[$.

3.2. *Transition rates as quality measure.* Let us now explain in detail our measure of quality designed to give a concept of optimal tuning which, as opposed to physical measures (see [13]), is robust for model reduction to Markov chains just retaining the jump dynamics between the equilibria of the diffusion. We shall use a notion that is based just on this rough transition mechanism. In fact, generalizing an approach for two state Markov chain models (see [8]), we measure the quality of tuning by computing for varying energy parameters $\mu$ the probability that, starting in $x_i$, the diffusion is transferred to $x_{-i}$ within the time window $[(a_\mu^i - h)T^\varepsilon, (a_\mu^i + h)T^\varepsilon]$ of width $2hT^\varepsilon$. To find the *stochastic resonance point* for large $T^\varepsilon$ (small $\varepsilon$), we have to maximize this measure of quality in $\mu \in I_R$. The probability for transition within this window will be approximated by the estimates of the preceding section. Uniformity of convergence to the exponential rates will enable us to maximize in $\mu$ for fixed small $\varepsilon$.

Let us now make these ideas precise. To make sure that the transition window makes sense at least for small $h$, we have to suppose that $a_\mu^i > 0$, $i = \pm$ for $\mu \in I_R$. This is guaranteed if

$$e_i(0) > \inf_{t \geq 0} \max_{i = \pm} e_i(t), \qquad i = \pm.$$

If this is not granted from the beginning, it suffices to start the diffusion a little later. For $\varrho$ small enough so that $B_\varrho(x_\pm) \subset A_\pm$, we call

$$(3.1) \quad \mathcal{M}(\varepsilon, \mu, \varrho) = \min_{i=\pm} \sup_{x \in B_\varrho(x_i)} \mathbb{P}_x(\tau_\varrho^{-i} \in [(a_\mu^i - h)T^\varepsilon, (a_\mu^i + h)T^\varepsilon]),$$

$$\varepsilon > 0, \mu \in I_R,$$

*transition probability for a time window of width $h$.* Here

$$\tau_\varrho^i = \inf\{t \geq 0 : X_t^\varepsilon \in B_\varrho(x_i)\}.$$

We are ready to state our main result on the asymptotic law of transition time windows.

THEOREM 3.2. *Let $M$ be a compact subset of $I_R$, $h_0 > 0$ and $\varrho$ be given according to Theorem* 2.3. *Then for all $h \leq h_0$,*

$$(3.2) \quad \lim_{\varepsilon \to 0} \varepsilon \log(1 - \mathcal{M}(\varepsilon, \mu, \varrho)) = \max_{i=\pm}\{\mu - e_i(a_\mu^i - h)\},$$

*uniformly for $\mu \in M$.*

PROOF. This is an obvious consequence of Theorem 2.3. □

It is clear that, for $h$ small, the eventually existing global minimizer $\mu_R(h)$ of

$$I_R \ni \mu \mapsto \max_{i=\pm}\{\mu - e_i(a_\mu^i - h)\}$$

is a good candidate for our resonance point. But it still depends on $h$. To get rid of this dependence, we shall consider the limit of $\mu_R(h)$ as $h \to 0$.

DEFINITION 3.3. *Suppose that*

$$I_R \ni \mu \mapsto \max_{i=\pm}\{\mu - e_i(a_\mu^i - h)\}$$

*possesses a global minimum $\mu_R(h)$. Suppose further that $\mu_R = \lim_{h \to 0} \mu_R(h)$ exists in $I_R$. We call $\mu_R$ the stochastic resonance point of the diffusion $(X^\varepsilon)$ with time periodic drift $b$.*

Let us now illustrate this resonance notion in a situation in which the energy functions are related by a phase lag $\phi \in ]0, 1[$, that is, $e_-(t) = e_+(t + \phi)$ for all $t \geq 0$. We shall show that in this case the stochastic resonance point exists if one of the energy functions, and thus both, has a unique point of maximal decrease on the interval where it is strictly decreasing (see Figure 4).

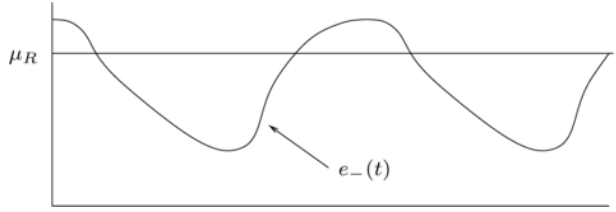

Fig. 4. *Point of maximal decrease.*

THEOREM 3.4. *Suppose that $e_-$ is twice continuously differentiable and has its global maximum at $t_1$, and its global minimum at $t_2$, where $t_1 < t_2$. Suppose further that there is a unique point $t_1 < s < t_2$ such that $e_-|_{]t_1,s[}$ is strictly concave, and $e_-|_{]s,t_2[}$ is strictly convex. Then $\mu_R = e_-(s)$ is the stochastic resonance point.*

PROOF. As a consequence of the phase lag of the energy functions,
$$\max_{i=\pm}\{\mu - e_i(a_\mu^i - h)\} = \{\mu - e_-(a_\mu^- - h)\}.$$

Write $a_\mu = a_\mu^-$ and recall that on the interval of decrease of $e_-$, $a_\mu = e_-^{-1}(\mu)$. The derivative of $\mu \mapsto \mu - e_-(a_\mu - h)$ has to vanish for the minimizer, which yields
$$1 = e'_-(a_\mu - h) \cdot a'_\mu = e'_-(a_\mu - h) \cdot \frac{1}{e'_-(a_\mu)}.$$

Our hypotheses concerning convexity and concavity of $e_-$ essentially means that $e''_-(s) = 0$, and $e''_-|_{]t_1,s[} < 0, e''_-|_{]s,t_2[} > 0$, which may be stated alternatively by saying that $\mu \mapsto e'_-(a_\mu)$ has a local maximum at $a_\mu = s$. Hence, for $h$ small, there exists a unique point $a_\mu(h)$ such that $e'_-(a_\mu(h) - h) = e'_-(a_\mu(h))$ and $\lim_{h\to 0} a_\mu(h) = s$. To show that $a_\mu(h)$ corresponds to a minimum of the function $\mu \mapsto [\mu - e_-(a_\mu - h)]$, we take the second derivative of this function at $a_\mu(h)$, which is given by
$$\frac{e'_-(a_\mu(h) - h)e''_-(a_\mu(h)) - e''_-(a_\mu(h) - h)e'_-(a_\mu(h))}{e'_-(a_\mu(h))}.$$

But $e'_-(a_\mu(h)), e'_-(a_\mu(h) - h) < 0$, whereas $e''_-(a_\mu(h) - h) > 0, e''_-(a_\mu(h)) < 0$. This clearly implies that $a_\mu(h)$ corresponds to a minimum of the function. But by definition, as $h \to 0$, $a_\mu(h) \to s$. Therefore, finally, $e_-(s)$ is the stochastic resonance point. □

3.3. *The robustness of stochastic resonance.* In the small noise limit $\varepsilon \to 0$, it seems reasonable to assume that the periodicity properties of the diffusion trajectories, caused by the periodic forcing the drift term exhibits, are essentially captured by a simpler, reduced stochastic process: a continuous time Markov chain which just jumps between two states representing the equilibria in the two domains of attraction. Jump rates correspond to the transition mechanism of the diffusion. This is just the reduction idea ubiquitous in the physics literature, and explained, for example, in [11]. We shall now show that in the small noise limit both models, diffusion and Markov chain, produce the same resonance picture, if quality of periodic tuning is measured by transition rates.

To describe the reduced model, let $e_\pm$ be the energy functions corresponding to transitions from $A_\mp$ to $A_\pm$ as before. Assume a phase locking of the two functions according to the previous subsection, that is, assume that $e_-(t) = e_+(t+\phi), t \geq 0$, with phase shift $\phi \in {]}0,1{[}$. So, let us consider a time-continuous Markov chain $\{Y_t^\varepsilon, t \geq 0\}$ taking values in the state space $S = \{-,+\}$ with initial data $Y_0^\varepsilon = -$. Suppose the infinitesimal generator is given by

$$G = \begin{pmatrix} -\varphi\left(\dfrac{t}{T^\varepsilon}\right) & \varphi\left(\dfrac{t}{T^\varepsilon}\right) \\ \psi\left(\dfrac{t}{T^\varepsilon}\right) & -\psi\left(\dfrac{t}{T^\varepsilon}\right) \end{pmatrix},$$

where $\psi(t) = \varphi(t+\phi), t \geq 0$, and $\varphi$ is a 1-periodic function describing a rate which just produces the transition dynamics of the diffusion between the equilibria $\pm$, that is,

$$(3.3) \qquad \varphi(t) = \exp\left\{-\frac{e_+(t)}{\varepsilon}\right\}, \qquad t \geq 0.$$

Note that by choice of $\varphi$,

$$(3.4) \qquad \psi(t) = \exp\left\{-\frac{e_-(t)}{\varepsilon}\right\}, \qquad t \geq 0.$$

Transition probabilities for the Markov chain thus defined are easily computed (see [8], Section 2). For example, the probability density of the first transition time $\sigma_i$ is given by

$$(3.5) \quad \begin{aligned} p(t) &= \varphi(t) \exp\left\{-\int_0^t \varphi(s)\,ds\right\}, & \text{if } i = -, \\ q(t) &= \varphi(t+\phi) \exp\left\{-\int_0^t \varphi(s+\phi)\,ds\right\}, & \text{if } i = +, \end{aligned}$$

$t \geq 0$. Equation (3.5) can be used to obtain results on exponential rates of the transition times $\sigma_i$ if starting from $-i$, $i = \pm$. We summarize them and apply them to the following measure of quality of periodic tuning:

$$(3.6) \quad \mathcal{N}(\varepsilon,\mu) = \min_{i=\pm} \mathbb{P}_i(\sigma_{-i} \in [(a_\mu^i - h)T^\varepsilon, (a_\mu^i + h)T^\varepsilon]), \qquad \varepsilon > 0, \mu \in I_R,$$

STOCHASTIC RESONANCE FOR MULTIDIMENSIONAL DIFFUSIONS   43

which is called *transition probability for a time window of width h* for the Markov chain.

Here is the asymptotic result obtained from a slight modification of Theorems 3 and 4 of [8] which consists of allowing more general energy functions than the sinusoidal ones used there and requires just the same proof.

THEOREM 3.5. *Let $M$ be a compact subset of $I_R$, $h_0 < \sup(a_\mu^{-1}, T/2 - a_\mu^{-1})$. Then for $0 < h \leq h_0$,*

$$(3.7) \qquad \lim_{\varepsilon \to 0} \varepsilon \ln(1 - \mathcal{N}(\varepsilon, \mu)) = \max_{i=\pm}\{\mu - e_-(a_\mu^i - h)\},$$

*uniformly for $\mu \in M$.*

It is clear from Theorem 3.5 that the reduced Markov chain $Y^\varepsilon$ and the diffusion process $X^\varepsilon$ have exactly the same resonance behavior. Of course, we may define the *stochastic resonance point* for $Y^\varepsilon$, just as we did for $X^\varepsilon$. So the following final robustness result holds true.

THEOREM 3.6. *The resonance points of $(X^\varepsilon)$ with time periodic drift $b$ and of $(Y^\varepsilon)$ with exponential transition rate functions $e_\pm$ coincide.*

**Acknowledgment.** We are much indebted to an anonymous referee for very constructive criticism.

S. HERRMANN
INSTITUT DE MATHÉMATIQUES ELIE CARTAN
UNIVERSITÉ HENRI POINCARÉ NANCY I
B.P. 239
54506 VANDOEUVRE-LÈS-NANCY CEDEX
FRANCE
E-MAIL: herrmann@iecn.u-nancy.fr

P. IMKELLER
D. PEITHMANN
INSTITUT FÜR MATHEMATIK
HUMBOLDT-UNIVERSITÄT ZU BERLIN
UNTER DEN LINDEN 6
10099 BERLIN
GERMANY
E-MAIL: imkeller@mathematik.hu-berlin.de
  peithman@mathematik.hu-berlin.de